\documentclass[draft,10pt,oneside,a4paper]{article} 		

\usepackage[fleqn]{amsmath} 							
\usepackage{amsthm,amsfonts,latexsym,amssymb,amscd} 	
\usepackage[mathscr]{eucal}   							
\usepackage[all]{xy}										

\usepackage[draft=false]{hyperref} 					

\newcommand{\B}{\mathcal{B}}
\newcommand{\F}{\mathcal{F}}
\newcommand{\G}{\mathcal{G}}
\renewcommand{\H}{\mathcal{H}} 	
\newcommand{\I}{\mathcal{I}}
\renewcommand{\L}{\mathcal{L}}		
\newcommand{\N}{\mathcal{N}}
\renewcommand{\O}{\mathcal{O}}	

\newcommand{\Eg}{\mathfrak{E}}\newcommand{\Fg}{\mathfrak{F}}
\newcommand{\Gg}{\mathfrak{G}}

\newcommand{\CC}{{\mathbb{C}}}

\newcommand{\NN}{{\mathbb{N}}}
\newcommand{\RR}{{\mathbb{R}}}
\newcommand{\TT}{{\mathbb{T}}}

\newcommand{\As}{{\mathscr{A}}}\newcommand{\Bs}{{\mathscr{B}}}\newcommand{\Cs}{{\mathscr{C}}}
\newcommand{\Ds}{{\mathscr{D}}}\newcommand{\Es}{{\mathscr{E}}}\newcommand{\Fs}{{\mathscr{F}}}
\newcommand{\Hs}{{\mathscr{H}}}\newcommand{\Is}{{\mathscr{I}}}

\newcommand{\Ms}{{\mathscr{M}}}\newcommand{\Ns}{{\mathscr{N}}}\newcommand{\Os}{{\mathscr{O}}}
\newcommand{\Ps}{{\mathscr{P}}}	
\newcommand{\Rs}{{\mathscr{R}}}

\newcommand{\Xs}{{\mathscr{X}}}\newcommand{\Ys}{{\mathscr{Y}}}

\DeclareFontFamily{U}{rsfs}{\skewchar\font127 }
\DeclareFontShape{U}{rsfs}{m}{n}{%
   <5> <6> rsfs5
   <7> rsfs7
   <8> <9> <10> <10.95> <12> <14.4> <17.28> <20.74> <24.88> rsfs10
}{}
\DeclareSymbolFont{rsfs}{U}{rsfs}{m}{n}
\DeclareSymbolFontAlphabet{\scr}{rsfs}

\newcommand{\Af}{\scr{A}}\newcommand{\Cf}{\scr{C}}\newcommand{\Df}{\scr{D}}\newcommand{\Ef}{\scr{E}}
\newcommand{\Mf}{\scr{M}}\newcommand{\Tf}{\scr{T}}

\DeclareMathOperator{\spa}{span}
\DeclareMathOperator{\id}{Id}
\DeclareMathOperator{\Ob}{Ob}

\DeclareMathOperator{\Hom}{Hom}

\DeclareMathOperator{\dom}{Dom}

\DeclareMathOperator{\Sp}{Sp}
\DeclareMathOperator{\ev}{ev}

\renewcommand{\emph}{\textbf} 						
\newcommand{\cj}[1]{\overline{#1}}					
\newcommand{\mip}[2]{(#1\mid #2)}					
\newcommand{\ip}[2]{\langle #1\mid #2\rangle}	
\newcommand{\imp}{\Rightarrow}						
\newcommand{\st}{\ : \ }									
\newcommand{\cs}{C*}

\newcommand{\hlink}[2]{\href{#1}{\texttt{#2}}} 

\newtheorem{theorem}{Theorem}[section]			

\newtheorem{proposition}[theorem]{Proposition}

\newtheorem{definition}[theorem]{Definition}

\numberwithin{equation}{section}  	
\title{\textbf{Non-Commutative Geometry, Categories and Quantum Physics}}

\author{\normalsize Paolo Bertozzini$^a$\footnote{Partially supported by the Thai Research Fund: grant n.~RSA4780022.},
Roberto Conti$^{b*}$\footnote{Current address: Dipartimento di Scienze, Universit\`a di Chieti-Pescara ``G. D'Annunzio'', Viale Pindaro 42, I-65127 Pescara, Italy.}, 
Wicharn Lewkeeratiyutkul$^{b*}$
\\
\normalsize  \textit{$^a$Department of Mathematics and Statistics, Faculty of Science and Technology}
\\
\normalsize \textit{Thammasat University, Bangkok 12121, Thailand}
\\
\normalsize e-mail: \texttt{paolo.th@gmail.com}
\\
\normalsize  \textit{$^{b}$Department of Mathematics and Computer Science}
\\
\normalsize \textit{Faculty of Science, Chulalongkorn University, Bangkok 10330, Thailand}
\\
\normalsize e-mail: \texttt{conti@sci.unich.it} 
\\
\normalsize e-mail: \texttt{Wicharn.L@chula.ac.th}
}

\date{\normalsize{08 November 2009, revised: 26 December 2011}}

\begin{document}

\maketitle

\begin{abstract}\noindent 
After an introduction to some basic issues in non-commutative geometry (Gel'fand duality, spectral triples), we present a ``panoramic view'' of the status of our current research program on the use of categorical methods in the setting of A.~Connes' 
non-commutative geometry: morphisms/categories of spectral triples, categorification of Gel'fand duality.
We conclude with a summary of the expected applications of ``categorical non-commutative geometry'' to structural questions in relativistic quantum physics: (hyper)covariance,  quantum space-time, (algebraic) quantum gravity.

\medskip

Quantum Physics, Space-Time. 

\medskip

\noindent
\emph{MSC-2000:}
					46L87,			
					46M15, 			
					16D90,			
					18F99, 			
					81R60,			
					81T05,			
 				  	83C65.			
\end{abstract}

\tableofcontents

\section{Introduction.}

The purpose of this review paper is to present the status of our research work on categorical non-commutative geometry and to contextualize it providing appropriate references.

The paper is organized as follows. In section~\ref{sec: categories} we introduce the basic elementary definitions about categories, functors, natural transformations and dualities just to fix our notation.

In section~\ref{sec: ncg}, we first provide a review of the basic dualities (Gelf'and, Serre-Swan and Takahashi) that constitute the main categorical motivation for  non-commutative geometry and then we pass to introduce the definition of A.~Connes spectral triple.

In the first part of section~\ref{sec: cncg}, we give an overview of our proposed definitions of morphisms between spectral triples and categories of spectral triples.
In the second part of section~\ref{sec: cncg} we show how to generalize Gel'fand duality to the setting of commutative full 
\hbox{\cs-categories} and we suggest how to apply this insight to the purpose of defining ``bivariant'' spectral triples as a correct notion of metric morphism.

The last section~\ref{sec: physics}, is mainly intended for an audience of mathematicians and tries to explain how categorical and non-commutative notions enter the context of quantum mathematical physics and how we hope to see such notions emerge in a non-perturbative treatment of quantum gravity.

The last part (section~\ref{sec: AQG}) is more speculative and contains a short overview of
our present research program in quantum gravity based on Tomita-Takesaki modular theory and categorical non-commutative geometry.

We have tried to provide an extensive biliography (updated till October 2009 and supplemented by a few additional references in appendix)  
in order to help to place our research in a broader landscape
and to suggest as much as possible future links with interesting ideas already developed.
Of course missing references are sole responsability of the ignorance of the authors, that are still trying to learn their way through the material. We will be grateful for any suggestion to improve the on-line version of the document.

\medskip

\emph{Notes and acknowledgments} The partial research support provided by
the Thai Research Fund (grant n.~RSA4780022) is kindly acknowledged.
The paper originates from notes prepared in occasion of a talk at the ``International Conference on Analysis and its Applications'' in Chulalongkorn University in May~2006. 
Most of the results have been announced in the form of research seminars
in Norway (University of Oslo), in Australia (ANU in Canberra, Macquarie University in Sydney, University of Queensland in Brisbane, La Trobe University in Melbourne, University of Newcastle) and in Italy (SISSA Trieste, Universit\`a di Roma II, Universit\`a di Bologna and Politecnico di Milano). One of the authors (P.B.) thanks Chulalongkorn University for the weekly hospitality during the last three years of research work.

\emph{Notes and acknowledgments for the revised version} 
A preliminar version of the paper appeared in the proceedings of the ``International Conference on Mathematics and Its Applications'' (ICMA-MU 2007) in Mahidol University in May 2007 and was subsequently published in a very shortened form in the special volume 2007 of East West Journal of Mathematics. 
The present paper is the second (and final) on-line version for the arXiv, updating and replacing the original submission in January 2008. 
It contains, apart from corrections of several typos, significant improvements in several sections: the bibliography has been updated to October 2009;  section~\ref{sec: physics} on applications to physics has been considerably expanded; references to some important developments (i.e.~those by A.~Connes on the reconstruction theorem and by B.~Mesland on ``$KK$-morphisms'' of spectral triples) have been added; an 
appendix at the end of the manuscript contains selected additional references appeared after October 2009. 

We thank Prof.~S.~J.~Summers and Prof.~W.~Lawton for reading the original manuscript and suggesting various improvements. 

\section{Categories.}\label{sec: categories}

Just for the purpose to fix our notation, we recall some general definitions on category theory, for a full introduction to the subject the reader can consult S.~MacLane~\cite{Mc} or M.~Barr-C.~Wells~\cite{BW}.

\subsection{Objects and Morphisms.}

A \emph{category} $\Cf$ consists of
\begin{itemize}
\item[a)]
a class\footnote{The family of objetcs can be a proper class. The category is called \emph{small} if the class of objects is actually a set.} of \emph{objects} $\Ob_\Cf$,
\item[b)]
for any two object $A,B\in \Ob_\Cf$ a set of \emph{morphisms} $\Hom_\Cf(A,B)$,
\item[c)]
for any three objects $A,B,C\in \Ob_\Cf$ a \emph{composition} map
\begin{equation*}
\circ: \Hom_\Cf(B,C)\times\Hom_\Cf(A,B)\to \Hom_\Cf(A,C)
\end{equation*}
that satisfies the following properties for all morphisms $f,g,h$ that can be composed:
\begin{gather*}
(f\circ g)\circ h=f\circ(g\circ h), 
\\
\forall A\in \Ob_\Cf, \ \exists \iota_A\in \Hom_\Cf(A,A) \st \iota_A\circ f=f, \ g\circ\iota_A=g.
\end{gather*}
\end{itemize}

A morphism $f\in\Hom_\Cf(A,B)$ is called an \emph{isomorphism} if there exists another morphism $g\in\Hom_\Cf(B,A)$ such that 
$f\circ g=\iota_B$ and $g\circ f=\iota_A$.

\subsection{Functors, Natural Transformations, Dualities.}

Given two categories $\Cf,\Df$, a \emph{covariant functor} $\Fg:\Cf\to\Df$ is a pair of maps
\begin{gather*}
\Fg:\Ob_\Cf\to\Ob_\Df, \quad \Fg: A\mapsto\Fg_A, \quad \forall A\in \Ob_\Cf, 
\\
\Fg:\Hom_\Cf\to\Hom_\Df,\quad \Fg: x\mapsto F(x), \quad \forall x\in \Hom_\Cf,
\end{gather*}
such that $x\in\Hom_\Cf(A,B)$ implies $\Fg(x)\in \Hom_\Df(\Fg_A,\Fg_B)$
and such that, for any two composable morphisms $f,g$ and any object $A$,
\begin{gather*}
\Fg(g\circ f)=\Fg(g)\circ \Fg(h), \qquad \Fg(\iota_A)=\iota_{\Fg_A}.
\end{gather*}

For the definition of a \emph{contravariant functor} we require $\Fg(x)\in \Hom_\Df(\Fg_B,\Fg_A)$,
whenever \hbox{$x\in \Hom_\Cf(A,B)$}.

A \emph{natural transformation} $\eta: \Fg\to \Gg$ between two functors
$\Fg,\Gg: \Cf\to \Df$, is a map
$\eta:\Ob_\Cf\to\Hom_\Df,\quad \eta: A\mapsto\eta_A\in\Hom_\Df(\Fg_A,\Gg_A)$,
such that the following diagram
\begin{equation*}
\xymatrix{
\Fg_A	\ar[r]^{\eta_A} \ar[d]_{\Fg(x)}	& 		\Gg_A	\ar[d]^{\Gg(x)}
\\
\Fg_B	\ar[r]_{\eta_B}							&		\Gg_B.
}
\end{equation*}
is commutative for all $x\in \Hom_\Cf(A,B)$, $A,B\in \Ob_\Cf$.
A natural transformation $\eta:\Fg\to\Gg$ is a \emph{natural isomorphism} (or natural equivalence) if $\eta_A$ is an isomorphism for all objects $A$; in this case we say that the functors $\Fg$ and $\Gg$ are naturally equivalent.

The functor $\Fg:\Cf\to \Df$ is
\begin{itemize}
\item
\emph{faithful} if, for all $A,B\in \Ob_\Cf$, its restriction to the set $\Hom_\Cf(A,B)$ is injective;
\item
\emph{full} if its restriction to $\Hom_\Cf(A,B)$ is surjective;
\item
\emph{representative} if for all $X\in \Ob_\Df$ there exists $A\in\Ob_\Cf$ such that $\Fg_A$ is isomorphic to $X$ in $\Df$.
\end{itemize}

A \emph{duality} (a contravariant equivalence) of two categories $\Cf$ and $\Df$ is a pair of contravariant functors $\Gamma: \Cf \to \Df$ and $\Sigma: \Df \to \Cf$ such that $\Gamma\circ\Sigma$ and $\Sigma\circ\Gamma$ are naturally equivalent to the respective identity functors $\I_\Df$ and $\I_\Cf$.
A duality is actually specified by two functors, but given any one of the two functors in the dual pair, the other one is unique up to natural isomorphism.
A functor $\Gamma$ is in a duality pair if and only if it is full, faithful and representative (see for example 
M.~Barr-C.~Wells~\cite[Definition~3.4.2]{BW}).
Categories that are in duality are considered ``essentially'' the same (modulo the reversing of arrows).

\medskip

Some important examples of ``geometrical categories'' i.e.~categories whose objects are sets equipped with a suitable structure, whose morphisms are ``structure preserving maps'' and with composition always given by the usual composition of functions are:
\begin{itemize}
\item 
sets and functions;
\item 
topological spaces and continuous maps;
\item 
differentiable manifolds and differentiable maps;
\item 
Riemannian manifolds (or also metric spaces) with global metric isometries;
\item 
Riemannian manifolds with Riemannian (totally geodesic) immersions/submersions;
\item 
orientable (Riemannian) $n$-dimensional manifolds with orientation preserving maps.\footnote{Note that, in general, it has no intrinsic meaning to say that a map between manifolds of different dimension preserve (or reverse) the orientation: a map between oriented manifolds, determines only a unique orientation for the normal bundle of the manifold.}
\end{itemize}

\begin{itemize}
\item[$\looparrowright$]
Problem: we are not aware of any definition in the literature of ``spin-preserving map'' between spin-manifolds of different dimension. In the case of manifolds with the same dimension, it is of course possible to say that a map preserves the spin-structure if there is an isomorphism (usually non-unique), between the pull-back of the spin-bundle of the target manifold and the spin-bundle on the source manifold, that ``intertwines'' the charge conjugation operators.
Anyway, even in this case, since spin-bundles are not ``natural bundles'' on a manifold,  there is no intrinsic notion of ``pull-back'' for spinor fields (unless we consider some special classes of manifolds such as K\"ahler spin-manifolds of a given dimension\footnote{P.~Bertozzini, R.~Conti, W.~Lewkeeratiyutkul,
Non-commutative (Totally Geodesic) Submanifolds and Quotient Manifolds, in preparation.}).

The correct solution of this problem (as in the case of ``orientation preserving'' maps) consists of equipping the morphisms (considered as ``relation submanifolds'' of the Cartesian product of the source and target (oriented) spin-manifolds) with their own additional ``spin-structure'' (orientation). Work on this issue is in progress~\footnote{P.~Bertozzini, R.~Conti, W.~Lewkeeratiyutkul, Categories of Spectral Triples and Morita Equivalence, work in progress.}.
\end{itemize}

Other examples of immediate interest for us include
\emph{vector bundles and bundle maps}, with composition of bundle maps
and
\emph{Hermitian vector bundles and (co)isometric bundle maps}.
For example, note that $K$-theory is the study of some special functors from the category of vector bundles to the category of
(Abelian) groups.

\section{Non-commutative Geometry (Objects).}\label{sec: ncg}

For an introduction to the subject we refer the readers to the books by A.~Connes~\cite{C}, G.~Landi~\cite{Lan}, H.~Figueroa-J.~Gracia-Bondia-J.~Varilly~\cite{FGV} (see also \cite{Var}) and M.~Khalkhali~\cite{Kha}; for spectral triples and their relation to index theory we also suggest A.~Rennie's lectures notes~\cite{Re6}. 

Non-commutative geometry, created by A.~Connes, is a powerful extension of the ideas of R.~Decartes' analytic geometry: to substitute ``geometrical objects'' with their Abelian algebras of functions; to ``translate'' the geometrical properties of spaces into algebraic properties of the associated algebras\footnote{A line of thought already present in J.L.~Koszul algebraization of differential geometry.} and to ``reconstruct'' the original geometric spaces as derived entities (the spectra of the algebras), a technique that appeared for the first time in the work of I.~Gel'fand on Abelian \cs-algebras 
in 1939.\footnote{Although similar ideas, previously developed by D.~Hilbert, are well known and used also in P.~Cartier-A.~Grothendieck's definition of schemes in algebraic geometry.}

Whenever such ``codifications'' of geometry in algebraic terms still make sense if the Abelian condition is dropped,\footnote{Usually in the non-commutative case, there are several inequivalent generalizations of the same condition for Abelian algebras.} we can simply work with non-commutative algebras considered as ``duals'' of ``non-commutative spaces''.

The existence of dualities between categories of ``geometrical spaces'' and categories ``constructed from Abelian algebras'' is the starting point of any generalization of geometry to the non-commutative situation. Here are some examples.

	\subsection{Non-commutative Topology.}

			\subsubsection{Gel'fand Theorem.}

For the details on operator algebras, the reader may refer to R.~Kadison-J.~Ringrose~\cite{KR}, M.~Takesaki~\cite{T} and B.~Blackadar~\cite{Bl}.
A complex unital \emph{algebra} $\As$ is a vector space over $\CC$ with an associative unital bilinear multiplication. $\As$ is \emph{Abelian} (commutative) if $ab=ba$, for all
$a,b\in \As$. An \emph{involution} on $\As$ is a conjugate linear map
$*: \As\to \As$ such that $(a^*)^*=a$ and $(ab)^*=b^*a^*$, for all $a,b\in \As$.
An involutive complex unital algebra is $\As$ called a \emph{\cs-algebra} if $\As$ is a Banach space with a norm $a\mapsto\|a\|$ such that $\|ab\|\leq\|a\|\cdot\|b\|$ and  $\|a^*a\|=\|a\|^2$, for all $a,b\in\As$.
Notable examples are the algebras of continuous complex valued functions $C(X;\CC)$ on a compact topological space with the 
``sup-norm'' and the algebras of linear bounded operators $\B(H)$ on the Hilbert space $H$.

\begin{theorem}[Gel'fand]\footnote{See for example~\cite[Theorems~II.2.2.4, II.2.2.6]{Bl} or~\cite[Section~6]{La0}}
\label{th: gel}
There exists a duality $(\Gamma^{(1)},\Sigma^{(1)})$ between the category $\Tf^{(1)}$, of continuous maps between compact Hausdorff topological spaces, and the category $\Af^{(1)}$, of unital homomorphisms of commutative unital \cs-algebras.
\end{theorem}
$\Gamma^{(1)}$ is the functor that associates to compact Hausdorff topological spaces 
\hbox{$X\in \Ob_{\Tf^{(1)}}$} the unital commutative \cs-algebras C$(X;\CC)$ of complex valued continuous functions on $X$ (with pointwise multiplication and conjugation and supremum-norm) and that to continuous maps $f: X\to Y$ associates the unital 
$*$-homomorphisms  \hbox{$f^\bullet: C(Y;\CC)\to C(X;\CC)$} given by the pull-back of continuous $\CC$-valued functions by $f$.

$\Sigma^{(1)}$ is the functor that associates to every unital commutative \cs-algebra $\As$ its spectrum
$\Sp (\As):=\{\omega \ | \ \omega: \As\to \CC \
\text{is a unital $*$-homomorphism}\}$
(as a topological space with the weak topology induced by the evaluation maps $\omega\mapsto\omega(x)$, for all $x\in \As$) and that to every unital $*$-homomorphism $\phi:\As\to \Bs$ of algebras associates the continuous map
$\phi^\bullet: \Sp(\Bs)\to \Sp(\As)$ given by the pull-back under $\phi$.

The natural isomorphism $\Gg:\I_{\Af^{(1)}}\to\Gamma^{(1)}\circ\Sigma^{(1)}$ is given by the \emph{Gel'fand transforms} $\Gg_\As: \As\to C(\Sp(\As))$ defined by
$\Gg_\As: a\mapsto \hat{a}$, where $\hat{a}:\Sp(\As)\to \CC$ is the Gel'fand transform of $a$
i.e.~$\hat{a}:\omega\mapsto \omega(a)$.

The natural isomorphism $\Eg: \I_{\Tf^{(1)}}\to\Sigma^{(1)}\circ\Gamma^{(1)}$ is given by the \emph{evaluation} homeomorphisms $\Eg_X: X\to \Sp(C(X))$ defined by
$\Eg_X: p\mapsto \ev_p$, where $\ev_p: C(X)\to\CC$ is the $p$-evaluation
i.e.~$\ev_p: f\mapsto f(p)$.

\smallskip

In view of this result, compact Hausdorff spaces and Abelian unital \cs-algebras are essentially the same thing and we can freely translate properties of the geometrical space
in algebraic properties of its Abelian algebra of functions.\footnote{For possible extensions of Gel'fand theorem to Tychonoff spaces and locally convex $*$-algebras see M.~Carri\'on-\'Alvarez~\cite{CA}. A Gel'fand duality theory for ordered topological spaces has been elaborated by F.~Besnard~\cite{Be2}.}

In the spirit of non-commutative geometry, we can simply consider non-Abelian unital \cs-algebras as ``duals'' of ``non-commutative compact Hausdorff topological spaces''.

			\subsubsection{Serre-Swan and Takahashi Theorems.}

A \emph{left pre-Hilbert-\cs-module} ${}_\As M$ over the unital \cs-algebra $\As$ (whose positive part is denoted by $\As_+:=\{x^*x \ | \ x\in \As\}$) is a unital left module $M$ over the unital ring $\As$ that is equipped with an $\As$-valued inner product $M\times M\to \As$ denoted by $(x,y)\mapsto {}_\As\ip{x}{y}$ such that, for all $x,y,z \in M$ and $a \in \As$,
$\ip{x+y}{z}=\ip{x}{z} + \ip{y}{z}$,
$\ip{a\cdot x}{z}=a\ip{x}{z}$,
$\ip{y}{x}=\ip{x}{y}^*$,
$\ip{x}{x}\in \As_+$,
$\ip{x}{x}=0_\As \imp x=0_M.$
A similar definition of a right pre-Hilbert-\cs-module is given with multiplication by elements of the algebra on the right.

A left Hilbert \cs-module ${}_\As M$ is a left pre-Hilbert \cs-module that is complete in the norm defined by $x\mapsto\sqrt{\|{}_\As\ip{x}{x}\|}$.\footnote{A similar definition applies for right modules.}
We say that a left pre-Hilbert \cs-module ${}_\As M$ is \emph{full} if
$\cj{\spa\{\ip{x}{y}\ | \ x,y\in M\}}=\As$, where the closure is in the norm topology of the \cs-algebra $\As$.
A \emph{pre-Hilbert-\cs-bimodule} ${}_\As M_\Bs$ over the unital \cs-algebras $\As,\Bs$, is a left pre-Hilbert module over $\As$ and a right pre-Hilbert \cs-module over $\Bs$ such that:
\begin{equation*}
(a\cdot x)\cdot b= a\cdot (x\cdot b), \quad \forall a\in \As, \ \forall x\in M, \ \forall b\in \Bs.
\end{equation*}
A full Hilbert \cs-bimodule is said to be an \emph{imprimitivity bimodule} or an \emph{equivalence bimodule} if:
\begin{equation*}
{}_\As\ip{x}{y}\cdot z=x\cdot \ip{y}{z}_\Bs, \quad \forall x,y,z\in M.
\end{equation*}
A bimodule ${}_\As M_\As$ is called \emph{symmetric} if $ax=xa$ for all $x\in M$ and
$a\in \As$.\footnote{Of course this definition make sense only for bimodules over a commutative algebra $\As$.}
A module ${}_\As M$ is \emph{free} if it is isomorphic to a module of the form $\oplus_J\As$ for some index set $J$. A module ${}_\As M$ is \emph{projective} if there exists another module ${}_\As N$ such that $M\oplus N$ is a free module.

\medskip

An ``equivalence result'' strictly related to Gel'fand theorem, is the following ``Hermitian'' version of Serre-Swan theorem (see for example M.~Karoubi~\cite[Theorem~6.18]{Kar} for the usual Serre-Swan equivalence and, for its Hermitian version, 
M.~Frank~\cite[Theorem~7.1]{Fr}, N.~Weaver~\cite[Theorem~9.1.6]{We2} and also  
H.~Figueroa-J.~Gracia-Bondia-J.~Varilly~\cite[Theorem~2.10 and page~68]{FGV}) that provides a ``spectral interpretation'' of symmetric finite projective Hilbert C*-bimodules over a commutative unital \cs-algebra as Hermitian vector bundles over the spectrum of the algebra.\footnote{\label{foot: 1} The result, as it is stated in the previously given references~\cite{Fr,We2} and~\cite[page~68]{FGV}, is actually formulated without the finitness and projectivity conditions on the modules and with Hilbert bundles (see J.~Fell-R.~Doran~\cite[Section~13]{FD} or~\cite[Definition~2.9]{FGV} for a detailed definition) in place of Hermitian bundles. Note that Hilbert bundles are not necessarily locally trivial, but they become so if they have finite constant rank (see for example J.~Fell-R.~Doran~\cite[Remark~13.9]{FD}) and hence the more general equivalence between the category of Hilbert bundles and the category of Hilbert \cs-modules actually entails the Hermitian version of Serre-Swan theorem presented here.}

\begin{theorem}[Serre-Swan]
Let $X$ be a compact Hausdorff topological space. Let $\Mf_{C(X)}$ be the category of symmetric projective finite Hilbert \cs-bimodules over the commutative \cs-algebra $C(X;\CC)$ with $C(X;\CC)$-bimodule morphisms.
Let $\Ef_X$ be the category of Hermitian vector bundles over $X$ with bundle morphisms\footnote{Continuous, fiberwise linear maps, preserving the base points.}.

The functor $\Gamma: \Ef_X\to\Mf_{C(X)}$, that to every Hermitian vector bundle associates
its symmetric $C(X)$-bimodule of sections, is an equivalence of categories.
\end{theorem}

In practice, to every Hermitian vector bundle $\pi: E\to X$ over the compact Hausdorff space $X$, we associate
the symmetric Hilbert \cs-bimodule $\Gamma(X;E)$, the continuous sections of $E$, over the \cs-algebra $C(X;\CC)$.

Since, in the light of Gel'fand theorem, non-Abelian unital \cs-algebras are to be interpreted as ``non-commutative compact Hausdorff topological spaces'', Serre-Swan theorem suggests that finite projective Hilbert \cs-bimodules over unital \cs-algebras should be considered as ``Hermitian bundles over non-commutative Hausdorff compact spaces''.
\begin{itemize}
\item[$\looparrowright$]
Problem: Serre-Swan theorem deals only with categories of bundles over a fixed topological space (categories of modules over a fixed algebra, respectively).
In order to extend the theorem to categories of bundles over different spaces, it is necessary to define generalized notions of morphism between modules over different algebras.
The easiest solution is to define a morphism from the $\As$-module ${}_\As \Ms$ to the 
$\Bs$-module ${}_\Bs \Ns$ as a pair $(\phi,\Phi)$, where $\phi: \As\to \Bs$ is a homomorphism of algebras and $\Phi:\Ms\to \Ns$ is a $\CC$-linear map of the bimodules such that $\Phi(a m)=\phi(a)\Phi(m)$, for all $a\in \As$ and $m\in \Ms$.
This is the notion that we have used in~\cite{BCL6} 
and that appeared also in~\cite{Ta1,Ta2,FGV,Ho}.
A more appropriate solution would be to consider ``congruences'' of bimodules and reformulate Serre-Swan theorem in terms of relators (as defined in~\cite{BCL6}). Work on this topic is in progress\footnote{P.~Bertozzini, R.~Conti, W.~Lewkeeratiyutkul,
Categories of Spectral Triples and Morita Equivalence, work in progress.}.
\item[$\looparrowright$]
Problem: note that Serre-Swan theorem gives an equivalence of categories (and not a duality), this will create problems of ``covariance'' for any  generalization of the well-known covariant functors between categories of manifolds and categories of their associated vector (tensor, Clifford) bundles, to the case of non-commutative spaces and their ``bundles''.
Again a more appropriate approach using relators should deal with this issue.
\end{itemize}
A first immediate solution to both the above problems is provided by Takahashi duality theorem below. Serre-Swan equivalence is actually a particular case of the following general (and surprisingly almost unnoticed) Gel'fand duality result that was obtained in 1971 by A.~Takahashi~\cite{Ta1,Ta2}.\footnote{Note that our Gel'fand duality result for commutative full \cs-categories (that we will present later in section~\ref{sec: hgel}) can be seen as ``strict''-$*$-monoidal version of Takahashi duality.}
In this formulation, one actually consider much more general \cs-modules and Hilbert bundles at the price of losing contact with $K$-theory; anyway (as described in the footonote~\ref{foot: 1} at page~\pageref{foot: 1}) the Hermitian version of Serre-Swan theorem can be recovered considering bundles with constant finite rank (over a fixed compact Hausdorff topological space).

\begin{theorem}[Takahashi]
There is a (weak $*$-monoidal) category ${}_\bullet\Mf$ of left Hilbert \cs-modules ${}_\As M, {}_\Bs N$ over unital commutative \cs-algebras, whose morphisms are given by pairs $(\phi,\Phi)$ where  $\phi:\As\to \Bs$ is a unital $*$-homomorphism of 
\cs-algebras and $\Phi: M\to N$ is a continuous additive map such that $\Phi(ax)=\phi(a)\Phi(x)$, for all $a\in \As$ and $x\in M$.

There is a (weak $*$-monoidal) category $\Ef$ of Hilbert bundles $(\Es,\pi,\Xs), (\Fs,\rho,\Ys)$ over compact Hausdorff topological spaces with morphisms given by pairs $(f,\F)$ with $f:\Xs\to \Ys$ continuous and $\F:f^\bullet(\Fs)\to \Es$ a continuous fiberwise linear map that satisfies $\pi\circ\F=\rho^f$, where $(f^\bullet(\Fs),\rho^f,\Xs)$ denotes the pull-back of the bundle $(\Fs,\rho,\Ys)$ under $f$.

There is a duality (of weak $*$-monoidal) categories given by the functor $\Gamma$ that associates to every Hilbert bundle $(\Es,\pi,\Xs)$ the set of sections $\Gamma(\Xs;\Es)$ and that to every morphism of bundles $(f,\F): (\Es,\pi,\Xs)\to(\Fs,\rho,\Ys)$ associates the morphism of modules $(f^\bullet,\Phi)$, where $\Phi$ is the map that to evey section $\sigma\in \Gamma(\Ys;\Fs)$ associates the section
$\F\circ f^\bullet(\sigma)\in \Gamma(\Xs;\Es)$.
\end{theorem}

Of course, much more deserves to be said about the vast landscape of research currently developing in non-commutative topology, but it is not our purpose to provide here an overview of this huge subject. Fairly detailed treatments of some of the usual techniques in algebraic topology are already available in their non-commutative counterpart (see~\cite{FGV} or the expository article by J.~Cuntz~\cite{Cu} for more details): non-commutative $K$-theory ($K$-theory of \cs-algebras), $K$-homology 
(G.~Kasparov's $KK$-theory) and (co)homology (Hochschild and A.~Connes, B.~Tsygan cyclic cohomologies). Among the most recent achievements, we limit ourselves to mention the extremely interesting definitions of quantum principal and associated bundles by P.~Baum-P.~Hajac-R.~Matthes-W.~Szymanski~\cite{BHMS} and of non-commutative CW-complexes by D.~N.~Diep~\cite{Di}.

At the (differential) topological level, we mention that important connections between non-commutative geometry and signal processing are emerging in the works by O.~Bratteli-P.~Jorgensen~\cite{BJ} (wavelets and Cuntz algebras) and by 
F.~Luef~\cite{Lu,Lu2} (Gabor analysis and Hilbert C*-modules for non-commutative tori).

	\subsection{Non-commutative (Spin) Differential Geometry.} 
	

What are ``non-commutative manifolds''?

In order to define ``non-commutative manifolds'', we have to find a categorical duality between a category of manifolds and a suitable category constructed out of Abelian \cs-algebras of functions over the manifolds.
The complete answer to the question is not yet known, but (at least in the case of compact  finite-dimensional orientable Riemannian spin-manifolds) the notion of Connes spectral triples and Connes-Rennie-Varilly~\cite{C4,C17,RV} reconstruction theorem provide an appropriate starting point, specifying the objects of our non-commutative category.\footnote{We will of course deal later with the morphisms in section~\ref{sec: msp}.}

			\subsubsection{Connes Spectral Triples.}
			
A.~Connes (see~\cite{C, FGV}) has proposed a set of axioms for ``non-commutative manifolds'' (at least in the case of a compact finite-dimensional orientable Riemannian spin-manifolds), called a (compact) spectral triple or an (unbounded) $K$-cycle.
\begin{itemize}
\item
A (compact) \emph{spectral triple} $(\As, \H, D)$ is given by:
\begin{itemize}
\item
a unital pre-\cs-algebra $\As$;\footnote{Sometimes $\As$ is required to be closed under holomorphic functional calculus.}
\item
a (faithful) representation $\pi: \As \to \B(\H)$ of $\As$ on the Hilbert space $\H$;
\item
a (generally unbounded) self-adjoint operator $D$ on $\H$, called the Dirac operator, such that:
\begin{itemize}
\item[a)]
the resolvent $(D-\lambda)^{-1}$ is a compact operator, $\forall \lambda \in \CC\setminus\RR$,\footnote{As already noticed by Connes,
this condition has to be weakened in the case of non-compact manifolds, cf.~\cite{GLMV,GGISV,Re3,Re4}.}
\item[b)]
$[D,\pi(a)]_{-}\in \B(\H)$, 
for every $a \in \As,$ \\
where $[x,y]_{-}:= xy- yx$ denotes the commutator of $x,y \in \B(\H)$.\footnote{Since the Dirac operator $D$ can be unbounded, the condition $[D,\pi(a)]_-\in \B(\H)$ actually means that the domain of $D$ is invariant under all the elements $a\in\pi(\As)$ and that the operators $[D,\pi(a)]_-=D\circ\pi(a)-\pi(a)\circ D$, defined on $\dom (D)\subset \H$, can be extended to bounded linear operators on $\H$.}
\end{itemize}
\end{itemize}
\item
A spectral triple is called \emph{even} if there exists a grading operator, i.e.~a bounded self-adjoint operator $\Gamma \in \B(\H)$ such that:
\begin{gather*}
\Gamma^2=\text{Id}_\H; \quad  [\Gamma, \pi(a)]_{-}=0, \ \forall a \in \As; \quad
[\Gamma, D]_{+}=0,
\end{gather*}
where $[x,y]_{+}:=xy+yx$ is the anticommutator of $x,y.$

A spectral triple that is not even is called \emph{odd}.
\item
A spectral triple is \emph{regular} if
the function
$$\Xi_x: t\mapsto \exp(it|D|)x\exp(-it|D|)$$
is regular,
i.e.~$\Xi_x\in \text{C}^\infty(\RR, \B(\H)),$\footnote{
This condition is equivalent to
$\pi(a), [D,\pi(a)]_{-} \in \cap_{m=1}^\infty \text{Dom}\, \delta^m,$ for all $a\in \As$, where $\delta$ is the derivation given by $\delta(x):=[|D|, x]_{-}$.}
for every $x \in \Omega_D(\As)$,
where~\footnote{We assume that for $n=0 \in \NN$ the term in the formula simply reduces to $\pi(a_0)$.}
\begin{equation*}
\Omega_D(\As):=\spa \{\pi(a_0)[D,\pi(a_1)]_- \cdots [D,\pi(a_n)]_- \  |
\ \  n\in \NN, \ a_0, \dots, a_n \in \As\} \ .
\end{equation*}
\item
A spectral triple is \emph{$n$-dimensional} iff there exists an integer $n$ such that the Dixmier trace of $|D|^{-n}$ is finite nonzero.
\item
A spectral triple is \emph{$\theta$-summable} if $\exp(-t D^2)$
is a trace-class operator for all $t>0$. 
\item
A spectral triple is \emph{real} if there exists an antiunitary operator $J: \H \to \H$ such that:
\begin{gather*}
[\pi(a), J\pi(b^*)J^{-1}]_{-}=0, \quad \forall a,b \in \As; \\
[\, [D, \pi(a)]_{-}, J\pi(b^*)J^{-1}]_{-}=0, \quad \forall a,b \in \As, \quad {\emph{first order condition};} \\
J^2=\pm\text{Id}_\H;  \quad [J,D]_{\pm}=0;
\quad \text{and, only in the even case,} \quad
[J,\Gamma]_{\pm}=0,
\end{gather*}
where the choice of $\pm$ in the last three formulas depends on the ``dimension'' $n$ of the spectral triple modulo $8$ in accordance to  the following table:
\begin{center}\label{tb: J}
\begin{tabular}{|l|c|c|c|c|c|c|c|c|}
\hline
$n$	&$0$	&$1$	&$2$	&$3$	&$4$	&$5$	&$6$	&$7$	\\
	\hline
$J^2=\pm\text{Id}_\H$	&	$+$	&	$+$		&	$-$		&	$-$		&	$-$		&	$-$		&	$+$		& $+$		 \\
	\hline
$[J,D]_{\pm}=0$				&	$-$	&	$+$		&	$-$		&	$-$		&	$-$		&	$+$		&	$-$		& $-$		 \\
\hline
$[J,\Gamma]_{\pm}=0$	&	$-$	&			&	$+$		&			&	$-$		&			&	$+$		& 		\\
\hline
\end{tabular}
\end{center}
\item
A spectral triple is \emph{finite} if
$\H_\infty := \cap_{k=1}^\infty \text{Dom}\, D^k$ is a finite projective $\As$-bimodule 
and \emph{absolutely continuous} if, there exists an Hermitian form $(\xi,\eta)\mapsto\mip{\xi}{\eta}$ on $\H_\infty$ such that, for all $a\in \As$, $\ip{\xi}{\pi(a)\eta}$ is the Dixmier trace of $\pi(a)\mip{\xi}{\eta}|D|^{-n}$.
\item
An $n$-dimensional spectral triple is said to be \emph{orientable} if there is a Hochschild cycle
$c=\sum_{j=1}^m a^{(j)}_0\otimes a^{(j)}_1\otimes \cdots \otimes a^{(j)}_n$ such that its ``representation'' on the Hilbert space $\H$, $\pi(c)=\sum_{j=1}^m\pi(a^{(j)}_0)[D, \pi(a^{(j)}_1)]_-\cdots[D,\pi(a^{(j)}_n)]_-$ is the grading operator in the even case or the identity operator in the odd case\footnote{In the following, in order to simplify the discussion, we will always refer to a ``grading operator'' $\Gamma$ that actually coincides with the grading operator in the even case and that is by definition the identity operator in the odd case.}.
\item
A real spectral triple is said to satisfy \emph{Poincar\'e duality} if its fundamental class in the $KR$-homology of 
$\As\otimes \As^{\rm op}$ induces (via Kasparov intersection product) an isomorphism between the $K$-theory $K_\bullet(\As)$ and the $K$-homology $K^\bullet (\As)$ of $\As$.\footnote{In \cite{RV} some of the axioms are reformulated in a different form, in particular this condition is replaced by the requirement that the C*-module completion of
$\H_\infty$ is a Morita equivalence bimodule between (the norm completions of) $\As$ and $\Omega_D(\As)$.}
\item
A spectral triple will be called \emph{Abelian} or commutative whenever $\As$ is Abelian.
\item
A spectral triple is \emph{irreducible} if there is no non-trivial closed subspace in $\H$ that is invariant for $\pi(\As), D, J, \, \Gamma$.
\end{itemize}

To every spectral triple $(\As,\H,D)$ there is a naturally associated quasi-metric\footnote{In general $d_D$ can take the value $+\infty$ unless the spectral triple is irreducible.} on the set of pure states $\Ps(\As)$, called Connes' distance and given for all pure states $\omega_1,\omega_2$ by:
\begin{equation*}
d_D(\omega_1,\omega_2):=\sup\{|\omega_1(x)-\omega_2(x)| \ | \ x \in \As, \ \|[D,\pi(x)]\|\leq 1\}.
\end{equation*}

\begin{theorem}[Connes; see e.g.~\cite{C,FGV}]\label{th: Co}
Given an orientable compact Riemannian spin $m$-dimensional differentiable manifold $M$, with a given complex spinor bundle $S(M)$, a given spinorial charge conjugation $C_M$ and a given volume form $\mu_M$,\footnote{Remember that an orientable manifolds admits two different orientations and that, on a Riemannian manifold, the choice of an orientation canonically determines a volume form $\mu_M$. Recall also~\cite{S} that a spin-manifold $M$ admits several inequivalent spinor bundles and for every choice of a complex spinor bundle $S(M)$ (whose isomorphism class define the spin$^c$-structure of $M$) there are inequivalent choices of spinorial charge conjugations $C_M$ that define, up to bundle isomorphisms, the spin-structure of $M$.} define:
\begin{itemize}
\item[]
$\As_M:=C^\infty(M;\CC)$ the algebra of complex valued regular functions on the differentiable manifold $M$,
\item[]
$\H_M:=$L$^2(M;S(M))$ the Hilbert space of ``square integrable'' sections of the given spinor bundle $S(M)$ of the manifold $M$ i.e.~the completion of the space $\Gamma^\infty(M;S(M))$ of smooth sections of the spinor bundle $S(M)$ equipped with the inner product given by $\ip{\sigma}{\tau}:=\int_M\ip{\sigma(p)}{\tau(p)}_p\,\text{d}\mu_M$, where $\ip{\,}{\,}_p$, with $p\in M$, is the unique inner product on $S_p(M)$ compatible with the Clifford action and the Clifford product.
\item[]
$D_M$ the Atiyah-Singer Dirac operator i.e.~the closure of the operator that is obtained by ``contracting'' the unique spinorial covariant derivative $\nabla^{S(M)}$ (induced on $\Gamma^\infty(M;S(M))$ by the Levi-Civita covariant derivative of $M$, see~\cite[Theorem~9.8]{FGV}) with the Clifford multiplication;
\item[]
$J_M$ the unique antilinear unitary extension $J_M:\H_M\to \H_M$ of the operator determined by the spinorial charge conjugation 
$C_M$ as $(J_M\sigma)(p):=C_M(\sigma(p))$ for
$\sigma\in \Gamma^\infty(M;S(M))$ and $p\in M$;
\item[]
$\Gamma_M$ the unique unitary extension on $\H_M$ of the operator given by fiberwise grading on $S_p(M)$, with $p\in M$.\footnote{The grading is actually the identity in odd dimension.}
\end{itemize}
The data $(\As_M,\H_M,D_M)$ define a spectral triple that is Abelian regular finite absolutely continuous $m$-dimensional real, with real structure $J_M$, orientable, with grading $\Gamma_M$, and that satisfies Poincar\'e duality.
\end{theorem}


\begin{theorem}[Connes~\cite{C4,C17}]
\label{th: Co-Re}
Let $(\As,\H,D)$ be an irreducible commutative real (with real structure $J$ and grading $\Gamma$) strongly regular\footnote{In the sense of~\cite[Definition~6.1]{C17}.} $m$-dimensional finite absolutely continuous orientable spectral triple satisfying Poincar\'e duality.
The spectrum of (the norm closure of) $\As$
can be endowed,
essentially in a unique way, with the structure of an $m$-dimensional
connected compact spin Riemannian manifold $M$ with an
irreducible complex spinor bundle $S(M)$, a charge conjugation $J_M$ and a grading $\Gamma_M$ such that:
$\As\simeq C^\infty(M; \CC)$, $\H\simeq\text{L}^2(M,S(M))$, $D \simeq D_M$, $J \simeq J_M$, $\Gamma \simeq \Gamma_M$.
\end{theorem}

\begin{itemize}
\item[$\looparrowright$]
A.~Connes first proved the previous theorem under the additional condition that $\As$ is already given as the algebra of smooth complex-valued functions over a differentiable manifold $M$, namely $\As = C^\infty(M;\CC)$, and conjectured~\cite[Theorem~6, Remark~(a)]{C10} \cite{C4} the result for general commutative pre-\cs-algebras $\As$.
A tentative proof of this last fact has been published by A.~Rennie~\cite{Re1}; some gaps were pointed out in the original argument, 
a different revised, but still incorrect, proof appears in~\cite{RV} (see also~\cite{RV2}) under some additional technical conditions.
Recently A.~Connes~\cite{C17} (see also~\cite{C12}) finally provided the missing steps in the proof of the result. 
\end{itemize}

As a consequence, there exists a one-to-one correspondence between unitary equivalence classes of spectral triples and connected compact oriented Riemannian spin-manifolds up to spin-preserving isometric diffeomorphisms.

Similar results are also available for spin$^c$-manifolds~\cite[Theorem~6, Remark~(e)]{C10}.

\subsection{Examples.}

Of course, the most inspiring examples of spectral triples (starting from those arising from Riemannian spin-manifolds) are contained in A.~Connes' book~\cite{C} and an updated account of most of the available constructions is contained in A.~Connes-M.~Marcolli's lecture notes~\cite{CM5}.
Here below we provide a short guide to some of the relevant literature:
\begin{itemize}
\item
Abelian spectral triples arising from the Atiyah-Singer Dirac operator on Riemannian spin-manifolds,  A.~Connes~\cite{C}, and classical compact homogeneous spaces, ~M.~Rieffel~\cite{Ri5}.
\item
Spectral triples for the non-commutative tori, A.~Connes~\cite{C}.
\item
Discrete spectral triples,
T.~Krajewski~\cite{Kr}, M.~Paschke-A.~Sitarz~\cite{PS}.
\item
Spectral triples from Moyal planes (these are examples of ``non-compact'' triples), V.~Gayral-J.M.~Gracia-Bondia-B.~Iochum-T.~Sch\"uker-J.~Varilly~\cite{GGISV}.
\item
Examples of Non-commutative Lorentzian Spectral Triples (following the definition given by A.~Strohmaier~\cite{Str}), W.~D.~van~Suijlekom~\cite{Sui}.
\item
Spectral Triples related to the Kronecker foliation
(following the general construction by A.~Connes-H.~Moscovici~\cite{CMo} of spectral triples associated to crossed product algebras related to foliations), R.~Matthes-O.~Richter-G.~Rudolph~\cite{MRR}.
\item
Dirac operators as multiplication by length functions on finitely generated discrete (amenable) groups, A.~Connes~\cite{C9}, M.~Rieffel~\cite{Ri}.
\item
$K$-cycles and (twisted) spectral triples arising from supersymmetric quantum field theories, A.~Jaffe-A.~Lesniewski-K.~Osterwalder~\cite{JLO,JLO2}, D.~Kastler~\cite{K},  A.~Connes~\cite{C}, D.~Goswami~\cite{Go2}; cyclic cocycles from super 
KMS-states in algebraic quantum field theory, D.~Buchholz-H.~Grundling~\cite{BG} and spectral triples on super-Virasoro algebras in conformal field theory, S.~Carpi-R.~Hillier-Y.~Kawahigashi-R.~Longo~\cite{CHKL}. 
\item
Spectral triples associated to quantum groups (in some case it is necessary to modify the first order condition involving the Dirac operator, requiring it to hold only up to compact operators),
P.~Chakraborty-A.~Pal~\cite{ChP, ChP2, ChP3, ChP4,ChP5,ChP6,ChP7,ChP8,ChP9}, D.~Goswami~\cite{Go1}, A.~Connes~\cite{C8},
L.~Dabrowski-G.~Landi-A.~Sitarz-W.~van~Suijlekom-J.~Varilly~\cite{DLSSV1,DLSSV2},
J.~Kustermans-G.~Murphy-L.~Tuset \cite{KMT}, S.~Neshveyev-L.~Tuset \cite{NT}; and also spectral triples associated to homogeneus spaces of quantum groups:
L.~Dabrowski~\cite{Da}, L.~Dabrowski-G.~Landi-M.~Paschke-A.~Sitarz~\cite{DLPS},
F.~D'Andrea-L.~Dabrowski~\cite{DD,DD2}, F.~D'Andrea-G.~Landi~\cite{DAL},  
F.~D'Andrea-L.~Dabrowski-G.~Landi~\cite{DDL,DDL2},
\cite{D} (the latter is ``twisted'' according to A.~Connes-H.~Moscovici~\cite{CMo2,Mos}).
\item
Non-commutative manifolds and instantons, A.~Connes-G.~Landi~\cite{CL},
L.~Dabrowski G.~Landi-T.~Masuda~\cite{DLM}, L.~Dabrowski-G.~Landi~\cite{DL},
G.~Landi~\cite{Lan3,Lan4}, G.~Landi-W.~van~Suijlekom~\cite{LS1,LS2}.
\item
Non-commutative spherical manifolds A.~Connes-M.~Dubois-Violette~\cite{CDV,CDV2,CDV3}.
\item
Spectral triples for some classes of fractal spaces, A.~Connes~\cite{C},
 \\
D.~Guido-T.~Isola~\cite{GI,GI2,GI3}, C.~Antonescu-E.~Christensen~\cite{AC},
E.~Christensen C.~Ivan-M.~Lapidus~\cite{CIL}.
\item
Spectral Triples for AF \cs-algebras, C.~Antonescu-E.~Christensen~\cite{AC}.
\item
Spectral triples in number theory: A.~Connes~\cite{C}, A.~Connes-M.~Marcolli~\cite{CM5},
R.~Meyer~\cite{Me};
spectral triples from Arakelov Geometry, from Mumford curves and hyperbolic Riemann surfaces,
C.~Consani-M.~Marcolli~\cite{CoM,CoM2,CoM3,CoM4},  G.~Cornelissen-M.~Marcolli-K.~Reihani-A.~Vdovina~\cite{CMRV}, G.~Cornelissen-M.~Marcolli~\cite{CMa}; spectral triples for certain classes of finite connected unoriented graphs, 
J.~W.~de Jong~\cite{DJ}. 
\item
Spectral triples of the standard model in particle physics, A.~Connes-J.~Lott~\cite{CLo},
J.~Gracia-Bondia-J.Varilly~\cite{GV}, D.~Kastler~\cite{K3,K4,KaS},
A.~Connes~\cite{C3,C4,C16},
J.~Barrett~\cite{Bar}, A.~Chamseddine-A.~Connes~\cite{CC,CC2,CC3,CC4}, A.~Chamseddine~\cite{Ch}, 
A.~Connes-M.~Marcolli~\cite{CM5,CM6},
A.~Chamseddine-A.~Connes-M.~Marcolli~\cite{CCMa}.
\end{itemize}

	\subsection{Other Spectral Geometries.}

In the last few years several others variants and extensions of ``spectral geometries'' have been considered or proposed:
\begin{itemize}
\item
Lorentzian spectral geometries: A.~Strohmaier~\cite{Str},
M.~Paschke-R.~Verch~\cite{PV2}, M.~Paschke-A.~Sitarz~\cite{PS2} and also M.~Borris-R.~Verch~\cite{BV},
\item
Riemannian non-spin: S.~Lord~\cite{Lo},
\item
Laplacian, K\"ahler: J.~Fr\"ohlich-O.~Grandjean-A.~Recknagel~\cite{FGR,FGR2,FGR3,FGR4} (for a study of non-commutative Laplace operators and elliptic partial differential equations in non-commutative geometry see J.~Rosenberg~\cite{Ros}), 
\item
Following works by M.~Breuer~\cite{Br1,Br2} on Fredholm modules on von Neumann algebras, M-T.~Benameur-T.~Fack~\cite{BF,BF2} and more recently in a remarkable series of papers
~\cite{CP,CPS1,CPS2,CPRS1,CPRS2,CPRS3,CPRS4,CRSS,BCPRSW,PaR,CPR,CPR2,CPR3,CPR4,CRT},  
M-T.~Benameur-A.~Carey-D.~Pask-J.~Phillips-A.~Rennie-F.~Sukochev-K.~Tong-K.~Wojciechowski (see also J.~Kaad-R.~Nest-A.~Rennie~\cite{KNR} and A.~Carey-S.~Neshveyev-R.~Nest-A.~Rennie~\cite{CNNR}), have been trying to generalize the formalism of Connes spectral triples when the algebra of bounded operators on the Hilbert space of the triple is replaced by a more general von Neumann algebra that is either semifinite or that carries a periodic action of the modular group of a KMS-state.

Among examples of semifinite spectral triples a special mention deserve those constructed on algebras of holonomy loops in canonical quantum gravity by J.~Aastrup-J.~Grimstrup-R.~Nest~\cite{AGN1,AGN2,AGN3,AGN4} (see also section~\ref{sec: qgrav}). 
\end{itemize}

\begin{itemize}
\item[$\looparrowright$]
Although non-commutative differential geometry, following A.~Connes, has been mainly developed in the axiomatic framework of spectral triples, that essentially generalize the structures available for the Atiyah-Singer theory of first order differential elliptic operators of the Dirac type, it is very likely that suitable ``spectral geometries'' might be developed using operators of higher order (the Laplacian type being the first notable example). Since ``topological obstructions'' (such us non-orientability, non-spinoriality)
are expected to survive essentially unaltered in the transition from the commutative to the non commutative world, these ``higher-order non-commutative geometries'' will deal with more general situations compared to usual spectral triples.
In this direction we are developing\footnote{P.~Bertozzini, R.~Conti, W.~Lewkeeratiyutkul,
Second Order Non-commutative Geometry, work in progress.}
definitions in the hope to obtain  Connes Rennie-Varilly reconstruction theorems also in these cases.
\end{itemize}

\begin{itemize}
\item[$\looparrowright$]
Apart from the ``spectral approaches'' to non-commutative geometry, more or less directly inspired by A.~Connes spectral triples, there are other lines of development that are worth investigating and whose ``relation'' with spectral triples is not yet clear:
\begin{itemize}
\item
J.-L.~Sauvageot~\cite{Sa} and F.~Cipriani~\cite{CS} are developing a version of non-commutative geometry described by Hilbert \cs-bimodules associated to a semigroup of completely positive contractions, an approach that is directly related to the analysis of the properties of the heat-kernel of the Laplacian on Riemannian manifolds
(see N.~Berline-E.~Getzler-M.~Vergne~\cite{BGV});
\item
M.~Rieffel~\cite{Ri4}, and along similar lines N.~Weaver~\cite{We1,We2}, have developed a theory of non-commutative compact metric spaces based on Lipschitz algebras.
\item
Following an idea of G.~Parfionov-R.~Zapatrin~\cite{PZ}, V.~Moretti~\cite{Mo} has generalized Connes' distance formula (using the D'Alembert operator) to the case of Lorentzian globally hyperbolic manifolds and has developed an approach to Lorentzian non-commutative geometry based on \cs-algebras whose relations with Strohmaier's spectral triples is intriguing.
\item
In algebraic quantum field theory (see section~\ref{sec: nc-st}), S.~Doplicher-K.~Fredenhagen~J.~Roberts~\cite{DFR0,DFR} (and also S.~Doplicher~\cite{Do1,Do2,Do}) have developed a model of Poincar\'e covariant quantum spacetime.
\item
O.~Bratteli and collaborators~\cite{B,BR} and more recently M.~Madore~\cite{Ma} have been approaching the definition of non-commutative differential geometries through modules of derivations over the algebra of ``smooth functions''.
\item
Strictly related to the previous approach there is a formidable literature (see for example S.~Majid~\cite{Maj1,Maj2}) on non-commutative geometry based on ``quantum groups'' structures (Hopf algebras).
\item
Most of the physics literature use the term non-commutative geometry to indicate non-commutative spaces obtained by a quantum ``deformation'' of a classical commutative space.
\end{itemize}
\end{itemize}

\section{Categories in Non-Commutative Geometry.}\label{sec: cncg}

After the discussion of ``objects'' in non-commutative geometry, we now shift our attention to some very tentative definitions of morphism of non-commutative spaces and of categories of non-commutative spaces.

In the first subsection we present morphisms of ``spectral geometries''. We limit our discussion essentially to the case of morphisms of A.~Connes spectral triples, although we expect that similar notions might be developed also for other spectral geometries.

In the second subsection we describe some other extremely important categories of ``non-commutative spaces'' that arise, at the ``topological level'', from ``variations on the theme'' of Morita equivalence.
Finally we indicate some direction of future research.

	\subsection{Morphisms of Spectral Triples.}\label{sec: msp}

Having described A.~Connes spectral triples and somehow justified the fact that spectral triples are a possible definition for ``non-commutative'' compact finite-dimensional orientable Riemannian spin-manifolds, our next goal here is to discuss definitions of ``morphisms'' between spectral triples and to construct categories of spectral triples (for further details and an updated overview of this line of research see also the slides~\cite{B2}). 

Even for spectral triples, there are actually several possible notions of morphism, according to the amount of ``background structure'' of the manifold that we would like to see preserved:\footnote{And also depending on the kind of topological properties that we would like to ``attach'' to our morphisms: orientation, spinoriality, \dots }
\begin{itemize}
\item
the metric, globally (isometries),
\item
the metric, locally (totally geodesic maps, in the differentiable case),
\item
the Riemannian structure,
\item
the differentiable structure,
\end{itemize}

			\subsubsection{Totally Geodesic Spin-Morphisms.}
			
This is the notion of morphism of spectral triples that we proposed in~\cite{BCL6}.

Given two spectral triples $(\As_j,\H_j,D_j),$ with $j=1,2,$ a
\emph{morphism of spectral triples} is a pair
\begin{gather*}
(\As_1,\H_1,D_1)\xrightarrow{(\phi,\Phi)} (\As_2,\H_2,D_2),
\end{gather*}
where $\phi: \As_1\to \As_2$ is a $*$-morphism between the pre-\cs-algebras $\As_1,\As_2$ and $\Phi: \H_1\to \H_2$ is a bounded\footnote{It might be necessary to relax this condition and to consider also cases in which $\Phi$ is unbounded.} 
linear map in $\B(\H_1;\H_2)$ that ``intertwines'' the representations $\pi_1, \pi_2\circ\phi$ and the Dirac operators $D_1, D_2:$
\begin{gather}
\pi_2(\phi(x))\circ \Phi=\Phi\circ \pi_1(x), \quad \forall x \in \As_1, \notag \\
D_2\circ\Phi= \Phi \circ D_1, \label{eq: Phi}
\end{gather}
i.e.~such that the following diagrams commute for every $x\in \As_1:$
\begin{equation*}
\xymatrix{
{\H_1} \ar[d]_{D_1} \ar[r]^{\Phi}\ar@{}[dr]|{\circlearrowleft} & {\H_2}\ar[d]^{D_2} \\
{\H_1}\ar[r]^{\Phi} & {\H_2}
}
\quad
\xymatrix{
{\H_1} \ar[d]_{\pi_1(x)} \ar[r]^{\Phi}\ar@{}[dr]|{\circlearrowleft} & {\H_2}\ar[d]^{\pi_2\circ\phi(x)} \\
{\H_1}\ar[r]^{\Phi} & {\H_2}
}
\end{equation*}
Here the intertwining relation between the Dirac operators holds on the domain of $D_1$, since we suppose that $\Phi(\dom (D_1))\subset \dom (D_2)$.

It is possible (in the case of even and/or real spectral triples) to require also commutations between $\Phi$ and the grading operators and/or the real structures.
More specifically:
\begin{itemize}
\item[]
a \emph{morphism of real spectral triples} $(\As_j,\H_j,D_j,J_j)$, is a morphism of spectral triples, as above, such that $\Phi$ also ``intertwines'' the real structure operators $J_1, J_2$: $J_2\circ \Phi=\Phi \circ J_1$;
\item[]
a \emph{morphism of even spectral triples} $(\As_j,\H_j,D_j,\Gamma_j)$, with $j=1,2$,
is a morphism of spectral triples, as above, such that $\Phi$ also ``intertwines'' the grading operators $\Gamma_1,\Gamma_2$: $\Gamma_2\circ \Phi=\Phi \circ \Gamma_1$.
\end{itemize}

Clearly this definition of morphism contains as a special case the notion of (unitary) equivalence of spectral triples~\cite[pp.~485-486]{FGV} and implies quite a strong relationship between the spectra of the Dirac operators of the two spectral triples.

Loosely speaking, for $\phi$ epi and $\Phi$ coisometric (respectively mono and isometric), in the commutative case\footnote{The details are developed in: P.~Bertozzini, R.~Conti, W.~Lewkeeratiyutkul, Non-commutative (Totally Geodesic) Submanifolds and Quotient Manifolds,
in preparation.}, one expects such definition to become relevant only for maps that ``preserve the geodesic structures'' (totally geodesic immersions and respectively totally geodesic submersions). Note that (already in the commutative case) these maps might not necessarily be metric isometries: totally geodesic maps are local isometries but not always global isometries (but we do not have a counterexample yet).

Furthermore these morphisms depend, at least in some sense, on the spin
structures:\footnote{In the case of morphisms of even real spectral triples, the map should preserve in the strongest possible sense the spin and orientation structures of the manifolds (whatever this might mean).} this ``spinorial rigidity'' (at least in the case of morphisms of real even spectral triples) requires that such  morphisms between spectral triples of different dimensions might be possible only when the difference in dimension is a multiple of~8.

It might be interesting to examine alternative sets of conditions on the pairs $(\phi,\Phi)$ that allow for example to formalize the notion of ``immersion'' of a non-commutative manifold into another with arbitrary higher dimension, avoiding the requirements coming from the spinorial structures.
Some preliminary considerations along similar lines have been independently proposed by A.~Sitarz~\cite{Si} in his habilitation thesis. There it was suggested that the appropriate morphisms satisfy some ``graded intertwining relations'' with the relevant operators, indicating the possibility to formalize suitable sign rules depending on the involved dimensions (modulo 8). We plan to elaborate on this topic elsewhere\footnote{P.~Bertozzini, R.~Conti, W.~Lewkeeratiyutkul, Morphism of Spectral Triples and Spin Manifolds,
work in progress.}.

			\subsubsection{Metric Morphisms.}

In~\cite{BCL11} we introduce the following notion of metric morphism.
Given two spectral triples $(\As_j,\H_j,D_j)$, with $j=1,2$, denote by $\Ps(\As_j)$ the sets of pure states over
(the norm completion of)
$\As_j$.
A \emph{metric morphism} of spectral triples
\begin{equation*}
(\As_1,\H_1,D_1)\xrightarrow{\phi} (\As_2,\H_2,D_2)
\end{equation*}
is by definition a unital epimorphism\footnote{Note that if $\phi$ is an epimorphism, its pull-back $\phi^\bullet$ maps pure states into pure states.}
$\phi : \As_1\to\As_2$ of pre-\cs-algebras whose pull-back
$\phi^\bullet : \Ps(\As_2)\to\Ps(\As_1)$ is an isometry, i.e.
\begin{equation*}
d_{D_1}(\phi^\bullet(\omega_1),\phi^\bullet(\omega_2))=d_{D_2}(\omega_1,\omega_2), \quad \forall \omega_1,\omega_2\in \Ps(\As_2).
\end{equation*}

This notion of metric morphism is ``essentially blind'' to the spin structures of the non-commutative manifolds (that in this case appears only as a necessary complication\footnote{Since it is possible to define functional distances using also Laplacian operators, we expect this notion to continue to make sense once a suitable notion of ``Laplacian non-commutative manifold'' is developed.}).

			\subsubsection{Riemannian Morphisms.}

A less rigid notion of morphism of spectral triples\footnote{P.~Bertozzini, R.~Conti, W.~Lewkeeratiyutkul, Morphisms of Non-commutative Riemannian Manifolds, in preparation. See also the slides~\cite{B2}.} (a definition that, for unitary maps, was introduced by R.~Verch and M.~Paschke~\cite{PV}) consists of relaxing the ``intertwining'' condition~\eqref{eq: Phi} between $\Phi$ and the Dirac operators, imposing only ``intertwining relations'' with the commutators of Dirac operators with elements of the algebras. In more detail: given two spectral triples $(\As_j,\H_j,D_j),$ with $j=1,2,$ a
\emph{Riemannian morphism of spectral triples} is a pair
\begin{gather*}
(\As_1,\H_1,D_1)\xrightarrow{(\phi,\Phi)} (\As_2,\H_2,D_2),
\end{gather*}
where $\phi: \As_1\to \As_2$ is a $*$-morphism between the pre-\cs-algebras $\As_1,\As_2$ and $\Phi: \H_1\to \H_2$ is a bounded linear map in $\B(\H_1;\H_2)$ that ``intertwines'' the representations $\pi_1, \pi_2\circ\phi$ and the commutators of the Dirac operators
$D_1, D_2$ with the elements $x\in \As_1,\phi(x)\in \As_2$:
\begin{gather*}
\pi_2(\phi(x))\circ \Phi=\Phi\circ \pi_1(x), \quad \forall x \in \As_1, \\
[D_2,\phi(x)]\circ\Phi= \Phi \circ [D_1,x], \quad \forall x\in \As_1,
\end{gather*}
i.e.~such that the following diagrams commute for every $x\in \As_1$:
\begin{equation*}
\xymatrix{
{\H_1} \ar[d]_{[D_1,x]} \ar[r]^{\Phi}\ar@{}[dr]|{\circlearrowleft} & {\H_2}\ar[d]^{[D_2,\phi(x)]} \\
{\H_1}\ar[r]^{\Phi} & {\H_2}
}
\quad
\xymatrix{
{\H_1} \ar[d]_{\pi_1(x)} \ar[r]^{\Phi}\ar@{}[dr]|{\circlearrowleft} & {\H_2}\ar[d]^{\pi_2\circ\phi(x)} \\
{\H_1}\ar[r]^{\Phi} & {\H_2}
}
\end{equation*}

In the commutative case, when $\phi$ is epi and $\Phi$ is coisometric (respectively mono and isometric), this definition is expected to correspond to the Riemannian isometries (respectively coisometries) of compact finite-dimensional orientable Riemannian 
spin-manifolds.

\begin{itemize}
\item[$\looparrowright$]
These notions of morphism of spectral triples are only tentative and more examples need to be tested.
As pointed out by A.~Rennie, it is likely that the ``correct'' definition of morphism will evolve, but it will surely reflect the basic structure suggested here.
At the ``topological level'' pair of maps $(\phi,\Phi)$ that intertwine the actions of the algebras on the respective Hilbert spaces (but not the Dirac operators or their commutators), have recently been used by P.~Ivankov-N.~Ivankov~\cite{II} for the definition of finite covering (and fundamental group) of a spectral triple.
\item[$\looparrowright$]
The several notions of morphism of spectral triples described above are not as general as possible. In a wider perspective,\footnote{P.~Bertozzini, R.~Conti, W.~Lewkeeratiyutkul,
Categories of Spectral Triples and Morita Equivalence, work in progress.} 
a morphism of spectral triples $(\As_j,\H_j,D_j)$, where $j=1,2$, might be formalized as a ``suitable'' functor $\Fs: {}_{\As_2}\Mf\to {}_{\As_1}\Mf$, between the categories ${}_{\As_j}\Mf$ of $\As_j$-modules, having ``appropriate intertwining'' properties with the Dirac operators $D_j$.
Now, under some ``mild'' hypothesis, by Eilenberg-Gabriel-Watts theorem (see for example~\cite{Me2}), any such functor is given by ``tensorization'' by a bimodule. These bimodules, suitably equipped with spectral data (as in the case of spectral triples), provide the natural setting for a general theory of morphisms of non-commutative spaces (see~\cite{B2} for some concrete proposal). 

In this direction we mention the notion of ``spectral correspondences'' developed by A.~Connes-M.~Marcolli~\cite{CM6} and further utilized in M.~Marcolli-A.~Zainy~\cite{MZ}. 
\end{itemize}

	\subsubsection{Morita Morphisms.}

In the previous subsections we described in some detail some proposed notions of morphism of ``non-commutative spaces'' (described as spectral triples) at the ``metric'' level.
A few other discussions of non-commutative geometry in a suitable categorical framework,
have already appeared in the literature in a more or less explicit form. Most of them deal essentially with morphisms at the ``topological level'' and are making use of the notion of Morita equivalence that we are going to introduce.
\begin{definition}
Two unital \cs-algebras $\As,\Bs$ are said to be \emph{strongly Morita equivalent} if there exists an imprimitivity bimodule ${}_\As X_\Bs$.
\end{definition}


It is a standard procedure in algebraic geometry, to define ``spaces'' dually by their ``spectra'' i.e.~by the categories of (equivalence classes of) representations of their algebras. Hence, for a given unital \cs-algebra $\As$, we consider its category
${}_\As \Mf$ of (isomorphism classes of) left \cs-Hilbert $\As$-modules with morphisms given by (equivalence classes of) $\As$-linear module maps.

Morphisms between these ``non-commutative spectra'' are given by covariant functors between the categories of modules.\footnote{This kind of ``ideology'' about categories of ``non-commutative spectra'' is very fashionable in ``non-commutative algebraic geometry'' (see for example M.~Kontsevich and A.~Rosenberg~\cite{KR1,KR2,R}).}

The Eilenberg-Gabriel-Watts theorem (see e.g.~\cite{Me2}) assures that under suitable conditions every functor $\Fg: {}_\As\Mf \to {}_\Bs\Mf$ coincides ``up to a natural equivalence'' with the functor given by left tensorization with a \cs-Hilbert
$\As$-$\Bs$-bimodule ${}_\Bs X_\As$ (with $X$ unique up to isomorphism of bimodules) i.e.:
\begin{equation*}
\Fg({}_\As E)\simeq{}_\Bs X_\As\otimes{}_\As E.
\end{equation*}

Y.~Manin~\cite{M} has been advocating the use of such ``Morita morphisms'' (tensorizations with Hilbert \cs-bimodules) as the natural notion of morphism of non-commutative spaces.
In~\cite{C3,C4,C5} A.~Connes already discussed how to transfer a given Dirac operator
using Morita equivalence bimodules and compatible connections on them,
thus leading to the concept of ``inner deformations'' of a spectral geometry underlying the ``transformation rule''
$\widetilde{D} = D + A + JAJ^{-1}$ (where $A$ denotes the ``connection'').
In our work\footnote{P.~Bertozzini, R.~Conti, W.~Lewkeeratiyutkul,
Categories of Spectral Triples and Morita Equivalence, work in progress.},
we try to define a strictly related category of spectral triples,
based on the notions of connection on a Morita morphism, that contains
``inner deformations'' as isomorphisms.

More specifically, given two spectral triples $(\As_j,\H_j,D_j),$ with $j=1,2,$
by a \emph{Morita-Connes} morphism of spectral triples, we mean a pair $(X,\nabla)$ where $X$ is Morita morphism from $\As_1$ to $\As_2$ i.e.~an $\As_2$-$\As_1$-bimodule that is a Hilbert \cs-module over $\As_2$ and $\nabla$ is a Riemannian connection on the bimodule $X$ (the Dirac operators are related to the connection $\nabla$ by the ``inner deformation'' formula). The composition 
$(X^3,\nabla^3)$ 
of two Morita-Connes morphisms $(X^1,\nabla^1)$ and $(X^2,\nabla^2)$ is defined by taking the tensor product 
$X^3:=X^1\otimes_{\As_2} X^2$ of the bimodules and taking the connection $\nabla^3$ on $X$ given by:
\begin{equation*}
\nabla^3(\xi_1\otimes\xi_2)(h_1):=\xi_1\otimes(\nabla^2\xi_2)(h_1)+(\nabla^1\xi_1)(\xi_2\otimes h_1), \quad h_1\in \H_1,\  \xi_j\in X^j. 
\end{equation*}

In a remarkable recent paper, A.~Connes-C.~Consani-M.~Marcolli~\cite{CCM} have been pushing even further the notion of ``Morita morphism'' defining morphisms between two algebras $\As,\Bs$ as ``homotopy classes'' of  bimodules in G.~Kasparov $KK$-theory $KK(\As,\Bs)$. In this way, every morphism is determined by a bimodule that is further equipped with additional structure (Fredholm module).\footnote{Other important results in this direction are obtained by S.~Mahanta~\cite{Mah4}.} 
In the same paper~\cite{CCM}, A.~Connes and collaborators provide ground for considering ``cyclic cohomology'' as an ``absolute cohomology of non-commutative motives'' and the category of modules over the ``cyclic category'' (already defined by 
A.~Connes-H.~Moscovici~\cite{CMo3}) as a ``non-commutative motivic cohomology''.

\begin{itemize}
\item[$\looparrowright$]
All the notions of categories of non-commutative spaces developed from the notion of Morita morphism, seem to be confined to the topological setting. Morita equivalence in itself is a non-commutative ``topological'' notion.
It is widely believed that Morita equivalent algebras should be considered as describing the ``same'' space. This comes from the fact that most of the ``geometric functors'' for commutative spaces when suitably extended to the non-commutative case are invariant under Morita equivalences (because Morita equivalence reduces to isomorphism for commutative algebras).
Anyway, most of the success of A.~Connes' non-commutative geometry actually comes from the fact that some commutative algebras are replaced with some other Morita equivalent non-commutative algebras that are able to describe in a much better way the geometry of the ``singular space''.

In a more direct way, it seems that the correct way to associate a \cs-algebra to a space, requires the direct input of the natural symmetries of the space (hence Morita equivalence is broken). Along these lines we have some work in progress on non-commutative Klein program\footnote{P.~Bertozzini, R.~Conti, W.~Lewkeeratiyutkul, Non-commutative Klein-Cartan Program, work in progress.}.

Although the formalization of the notion of morphism as a bimodule is probably here to stay, additional structures on the bimodule will be required to account for different level of ``rigidity'' (metric, Riemannian, differential, \dots) and some of these, are probably going to break Morita equivariance as long as non-topological properties are concerned.
\end{itemize}


A.~Connes-M.~Marcolli~\cite[Chapter~8.4]{CM6}  
and M.~Marcolli-A.~Zainy~\cite{MZ}  
give a definition of ``spectral correspondences'' as Hilbert C*-bimodules providing a ``bivariant version'' of a spectral triple. 

The problem of defining a ``metric'' category of spectral triples via morphisms in Kasparov $KK$-theory suitably equipped with smooth and metric structures, has been recently addressed in a remarkable paper by 
B.~Mesland~\cite{Mes}: 
a morphism from the  spectral triple $(\Bs,\Hs',D')$ to the spectral triple $(\As,\Hs,D)$ is given by a unitary isomorphism class of an unbounded ``smooth'' $\As$-$\Bs$-bimodule $(\Es,S,\nabla)$ with connection $\nabla$ such that:
\begin{itemize}
\item
$[\nabla,S]$ is a completely bounded operator, 
\item
$\Hs$ is isomorphic to $\Es\otimes_\Bs\Hs'$,  
\item
$D=S\otimes \id + \id\otimes_\nabla D'$ with $\id\otimes_\nabla D'(e\otimes f):=(-1)^{\partial e}(e\otimes D'f+\nabla_{D'}(e)f)$. 
\end{itemize}

S.~Mahanta~\cite{Mah4} 
is trying to relate ``spectral correspondences'' with the ``geometric morphisms'' of derived categories of the differential graded categories already used in the non-commutative algebraic geometry approach to non-commutative spaces~\cite{Mah1,Mah2,Mah3}. 


\begin{itemize} 
\item[$\looparrowright$]
Finally we note that we have not been discussing here the role of quantum groups as possible symmetries of spectral triples (see for example the recent papers by D.~Goswami~\cite{Go3,Go4,Go5,Go6} and 
J.~Bhowmick-D.~Goswami-A.~Skalski~\cite{BG1,BG2,BG3,BG4,BG5,BGS} discussing quantum isometries of spectral triples).
\end{itemize}

	\subsection{Categorification (Topological Level).}

Categorification is the term, introduced by L.~Crane-D.~Yetter~\cite{CY}, to denote the generic process to substitute ordinary algebraic structures with categorical counterparts. The term is now mostly used to denote a wide area of research 
(see J.~Baez-J.~Dolan~\cite{BD2}) whose purpose is to use higher order categories to define categorial analogs of algebraic structures. This \emph{vertical categorification}\footnote{In general a $n$-category get replaced with a $n+1$-category, increasing the ``depth'' of the available morphisms, hence the terminology ``vertical'' adopted here.} is usually done by promoting sets to categories, functions to functors, \dots hence replacing a category with a $2$-category and so on.
In non-commutative geometry, where usually spaces are defined ``dually'' by ``spectra'' i.e.~categories of representations of their algebras of functions, this is a kind of compulsory step: morphisms of non-commutative spaces are actually particular functors between ``spectra''. In this sense, non-commutative geometry (and also ordinary commutative algebraic geometry of schemes) is already a kind of vertical categorification.

There are also more ``trivial'' forms of \emph{horizontal categorification} in which ordinary algebraic unital associative structures are interpreted as categories with only one object and suitable analog categories with more than one object are defined.
In this case the passage is from endomorphisms of a single object to morphisms between different objects\footnote{Hence the name ``horizontal'', adopted here, that implies that no jump in the ``depth'' of morphisms is required. J.~Baez~\cite{B} prefers to use the term \emph{oidization} for this case.}:
\begin{center}
\begin{tabular}{|l|l|}
\hline
Monoids 	& Small Categories (Monoidoids)	\\
\hline
Groups 		& Groupoids 				\\
\hline
Associative Unital Rings & Ringoids \\
\hline
Associative Unital Algebras & Algebroids \\
\hline
Unital \cs-algebras & \cs-categories (\cs-algebroids)\\
\hline
\end{tabular}
\end{center}

It is an extremely interesting future topic of investigation to discuss the interplay between ideas of categorification and non-commutative geometry \dots here we are really only at the beginning of a long journey and we can present only a few ideas.\footnote{Other approaches to the abstract concept of ``categorification''
have turned out to be useful in the theory of knots and links, see \cite{Kh1,Kh2}.}

			\subsubsection{Horizontal Categorification of Gel'fand Duality.}\label{sec: hgel}

As a first step in the development of a ``categorical non-commutative geometry'', we have been looking at a possible ``horizontal categorification'' of Gel'fand duality (theorem~\ref{th: gel}).
In practice, the purpose is:
\begin{itemize}
\item
to find ``suitable embedding functors'' $F:\Tf^{(1)}\to \Tf$ and $G:\Af^{(1)}\to \Af$ of the categories $\Tf^{(1)}$ (of compact Hausdorff topological spaces) and $\Af^{(1)}$ (of unital commutative \cs-algebras) into two categories $\Tf$ and $\Af$;
\item
to extend the categorical duality $(\Gamma^{(1)},\Sigma^{(1)})$ between $\Tf^{(1)}$ and $\Af^{(1)}$ provided by Gel'fand theorem, to a categorical duality between $\Tf$ and $\Af$ in such a way that the following diagrams are commutative up to natural isomorphisms
$\eta, \xi$:
\end{itemize}
\begin{equation*}
\xymatrix{
\Tf^{(1)} \ar@{->}[d]_F \ar@^{->}[rr]^{\Gamma^{(1)}} & & \Af^{(1)} \ar@^{->}[ll]^{\Sigma^{(1)}} \ar@{->}[d]^G & &
F \circ \Sigma^{(1)} \ar[r]^\eta & \Sigma\circ G,
\\
\Tf \ar@^{->}[rr]^\Gamma & & \ar@^{->}[ll]^\Sigma \Af, & &
G\circ \Gamma^{(1)} \ar[r]_\xi & \Gamma\circ F.
}
\end{equation*}

Since $\Af^{(1)}$ is a full subcategory of the category of \cs-algebras, we identify the horizontal categorification of $\Af^{(1)}$ as a subcategory of the category of small \cs-categories.

In~\cite{BCL17}, in the setting of \cs-categories,
we provide a definition of ``spectrum'' of a commutative full \cs-category as a one-dimensional unital Fell bundle over a suitable groupoid (equivalence relation) and we prove a categorical Gel'fand duality theorem generalizing the usual Gel'fand duality between the categories of Abelian \cs-algebras and compact Hausdorff spaces.

\medskip

As a byproduct, in~\cite{BCL17/2} we also obtain the following spectral theorem for imprimitivity bimodules over Abelian unital \cs-algebras: every such bimodule is obtained by ``twisting'' (by the two projection homeomorphisms) the symmetric bimodule of sections of a unique Hermitian line bundle over the graph of a unique homeomorphism between the spectra of the two \cs-algebras.

\begin{theorem}(P.~Bertozzini-R.~Conti-W.~Lewkeeratiyutkul~\cite[Theorem~3.1]{BCL17/2})  
Given an imprimitivity Hilbert \cs-bimodule ${}_\As M_\Bs$ over the Abelian unital \cs-algebras $\As,\Bs$, there exists a canonical homeomorphism\footnote{$R_{BA}$ is a compact Hausdorff subspace of $\Sp(\As)\times\Sp(\Bs)$ homeomorphic to $\Sp(\As)$ (resp.~$\Sp(\Bs)$) via the  projections $\pi_A: R_{BA}\to \Sp(\As)$ (resp.~$\pi_B:R_{BA}\to\Sp(\Bs)$).}
$R_{BA}: \Sp(\As)\to\Sp(\Bs)$ and a Hermitian line bundle $E$ over $R_{BA}$ such that ${}_\As M_\Bs$ is isomorphic to the (left/right) ``twisting''\footnote{If $M$ is a left module over $\Cs$ and $\phi:\As\to\Cs$ is an isomorphism, the left twisting of $M$ by $\phi$ is the module over $\As$ defined by $a\cdot x:=\phi(a)x$ for $a\in \As$ and $x\in M$.} of the symmetric bimodule
$\Gamma(R_{BA}; E)_{C(R_{BA};\CC)}$ of sections of the bundle $E$ by the two ``pull-back'' isomorphisms $\pi_A^\bullet:\As\to C(R_{BA};\CC)$, $\pi_B^\bullet:\Bs\to C(R_{BA};\CC)$.
\end{theorem}

\begin{itemize}
\item[$\looparrowright$]
This reconstruction theorem for imprimitivity bimodules is actually only the starting point for the development of a complete ``bivariant'' version of Serre-Swan equivalence and Takahashi duality. In this case we will generalize the previous spectral theorem to (classes of) bimodules over commutative unital \cs-algebras that are more general than imprimitivity bimodules; furthermore the appropriate notion of morphism will be introduced in order to get a categorical duality. We plan to return to this subject elsewhere\footnote{P.~Bertozzini, R.~Conti, W.~Lewkeeratiyutkul, Bivariant Serre-Swan Equivalence, in preparation.}.
\end{itemize}

A \emph{\cs-category}~\cite{GLR,Mit} is a category $\Cs$ such that the sets $\Cs_{AB}:=\Hom_\Cs(B,A)$ are complex Banach spaces and the compositions are bilinear maps,
there is an involutive antilinear contravariant functor $*:\Hom_\Cs\to\Hom_\Cs$
acting identically on the objects such that
$x^*x$ is a positive element in the $*$-algebra $\Cs_{AA}$ for every $x\in \Cs_{BA}$
(that is, $x^* x = y^* y$ for some $y \in \Cs_{AA}$),
$\|xy\|\leq\|x\|\cdot \|y\|, \ \forall x\in \Cs_{AB}, \ y\in \Cs_{BC}$,
$\|x^*x\|=\|x\|^2, \ \forall x\in \Cs_{BA}$.

In a \cs-category $\Cs$, the sets $\Cs_{AA}:=\Hom_\Cs(A,A)$ are unital \cs-algebras for all $A\in \Ob_\Cs$.
The sets $\Cs_{AB}:=\Hom_\Cs(B,A)$ have a natural structure of unital Hilbert \cs-bimodule on the \cs-algebras $\Cs_{AA}$ on the right and $\Cs_{BB}$ on the left.

A \cs-category is \emph{commutative} if the \cs-algebras $\Cs_{AA}$ are Abelian for all $A\in \Ob_\Cs$.
The \cs-category $\Cs$ is \emph{full} if all the bimodules $\Cs_{AB}$ are
full\footnote{In this case $\Cs_{AB}$ are imprimitivity bimodules.}.
A basic example is the \cs-category of linear bounded maps between Hilbert spaces.

\medskip

A \emph{Banach bundle}~\cite[Section~I.13]{FD} $(E,p,X)$ is given by a continuous open surjection $p:E\to X$ of Hausdorff topological spaces, whose total space $E$ is equipped with a continuous partial operation of addition $+:\{(e_1,e_2) \ | \ p(e_1)=p(e_2)\}\to E$, a continuous operation of multiplication by scalars $\cdot: \CC\times E\to E$ and a continuous norm $\|\cdot\|:E\to \RR$, making all the fibers $E_x:=p^{-1}(x)$ Banach spaces and such that, for all $x\in X$, the sets of the form $B_{U,\epsilon}:=\{e\in E \ | \ p(e)\in U, \|e\|<\epsilon\}$, where $\epsilon>0$ and $U$ is a neighbourhood of $x\in X$, constitute a base of neighbourhoods of $0_x\in E_x$ in the topology of $E$.


If the topological base space $X$ is equipped with the algebraic structure of category (let $X^o$ be the set of its units, let $r,s : X \to X^o$ be its range and source maps and  let $X^n:=\{(x_1,\dots, x_n)\in \times_{j=1}^n X\ | \ s(x_j)=r(x_{j+1})\}$ be its set of 
$n$-composable morphisms), we further require that the composition $\circ:X^2\to X$ is a continuous map.

If $X$ is an \emph{involutive category} (also known as a $*$-category~\cite{GLR,Mit} or a ``dagger category''~\cite{Sel,AbC}) 
i.e.~there is a map $*:X\to X$ with the properties $(x^*)^*=x, \ \forall x\in X$ and $(x\circ y)^*=y^*\circ x^*$, for all $(x,y)\in X^2$, we also require $*$ to be continuous. An involutive category $X$ is called an \emph{involutive inverse category} if $x\circ x^*\circ x=x$ for all $x\in X$.

\medskip

A \emph{Fell bundle}\footnote{Fell bundles over topological groups were first introduced by J.~Fell~\cite[Section~II.16]{FD} and later generalized to the case of groupoids by S.~Yamagami (see A.~Kumjian~\cite{Ku} or P.~Muhly-D.~Williams~\cite{MW} and references therein) and to the case of inverse semigroups by N.~Sieben (see R.~Exel~\cite[Section~2]{Ex}).} 
\emph{over the involutive inverse category} $X$ (see also~\cite{BCL17}) is a Banach bundle $(E,p,X)$ whose total space $E$ is equipped with a multiplication defined on the set $E^2:=\{(e,f)\ | \ (p(e),p(f))\in X^2\}$, denoted by $(e,f)\mapsto ef$,
and an involution $*:E\to E$ such that
\begin{gather*}
e(fg)=(ef)g,\quad \forall (p(e),p(f),p(g))\in X^3, \\
p(ef)=p(e)\circ p(f), \quad \forall e,f\in E^2, \\
\forall x,y \in X^2, \ \text{the restriction of $(e,f)\mapsto ef$ to
$E_x\times E_y$ is bilinear}, \\
\|ef\|\leq \|e\|\cdot\|f\|, \quad \forall e,f\in E^2, \\
(e^*)^*=e,\quad \forall e\in E, \\
p(e^*)=p(e)^*, \quad \forall e\in E, \\
\forall x\in X, \ \text{the restriction of $e\mapsto e^*$ to $E_x$ is conjugate linear}, \\
(ef)^*=f^*e^*,\quad \forall e,f\in E^2, \\
\|e^*e\|=\|e\|^2, \quad \forall e\in E, \\                          
e^*e\geq 0, \quad \forall e\in E,                                    
\end{gather*}
where, in the last line we mean that $e^*e$ is a positive element in the \cs-algebra $E_{p(e^*e)}$. 

It is in fact easy to see that for every $x\in  X^o$, and more generally for every Hermitian idempotent $x=x\circ x=x^*\in X$,  the fiber $E_x$ is a \cs-algebra.
A Fell bundle $(E,p,X)$ is said to be \emph{unital} if the \cs-algebras $E_x$, for
$x\in X^o$, are unital.
Note that the fiber $E_x$ has a natural structure of Hilbert \cs-bimodule over the \cs-algebras $E_{r(x)}$ on the left and $E_{s(x)}$ on the right.
A Fell bundle is said to be \emph{saturated} if the above Hilbert \cs-bimodules $E_x$ are full.
Note also that in a saturated Fell bundle, the Hilbert \cs-bimodules $E_x$ are imprimitivity bimodules.

\medskip

Let $\O$ be a set and $X$ a compact Hausdorff topological space.
\\
We denote by
$\Rs_\O:=\{(A,B) \ | \ A,B\in \Os\}$
the ``total'' equivalence relation in $\Os$ and by
$\Delta_X:=\{(p,p)\ | \ p\in X\}$
the ``diagonal'' equivalence relation in $X$.
\begin{definition}
A \emph{topological spaceoid}\footnote{Note that, despite the name and the involvement of groupoids, spaceoids are not directly related with the fractaloids introducted by I.~Cho-P.~Jorgensen~\cite{CJ}: our spaceoids are groupoids but are equipped with a suitable bundle structure and fractaloids (graph groupoids with fractal properties) are not a horizontal categorification of self-similar fractal spaces.} 
$(\Es,\pi,\Xs)$ is a saturated unital rank-one Fell bundle over the product involutive topological category $\Xs:=\Delta_X\times\Rs_\O$.
\end{definition}

Let $(\Es_j,\pi_j,\Xs_j)$, for $j=1,2$, be two spaceoids (where $\Xs_j=\Delta_{X_j}\times \Rs_{\O_j}$, with $\O_j$ sets and $X_j$ compact Hausdorff topological spaces for $j=1,2$). 
\begin{definition}
A morphism of spaceoids $(\Es_1,\pi_1,\Xs_1)\xrightarrow{(f,\F)}(\Es_2,\pi_2,\Xs_2)$ is a pair $(f,\F)$ where
\begin{itemize}
\item
$f:=(f_\Delta,f_\Rs)$ with $f_\Delta:\Delta_1\to\Delta_2$ a continuous map of topological spaces and $f_\Rs:\Rs_1\to\Rs_2$ an isomorphism of equivalence relations;
\item
$\F:f^\bullet(\Es_2)\to\Es_1$ is a fiberwise linear continuous $*$-functor such that $\pi_1\circ\F=(\pi_2)^f$, where 
$(f^\bullet(\Es_2),\pi_2^f,\Xs_1)$ denotes a given choice of an $f$-pull-back\footnote{Here we denote by $\pi_2^f: f^\bullet(\Es_2)\to\Xs_1$ the projection of the pull-back bundle $(f^\bullet(\Es_2),\pi_2^f,\Xs_1)$ and by $f^{\pi_2}:f^\bullet(\Es_2)\to\Es_2$ the morphism of bundles such that $\pi_2\circ f^{\pi_2}=f\circ\pi^f$.} of $(\Es_2,\pi_2,\Xs_2)$.
\end{itemize}

\end{definition}

Topological spaceoids constitute a category if composition is defined by
\begin{equation*}
(g,\G)\circ(f,\F):=(g\circ f, \F\circ f^\bullet(\G)\circ \Theta), 
\end{equation*}
where $\Theta$ is the natural isomorphism from $f^\bullet(g^\bullet(\Es_3))$ to $(g\circ f)^\bullet(\Es_3)$, and (having chosen 
$(\Es,\pi,\Xs)$ to be the $\iota_\Xs$-pull-back of itself) with identities given by
\begin{equation*}
\iota(\Es,\pi,\Xs):=(\iota_{\Xs},\iota_{\Es}).
\end{equation*}


\medskip

The category $\Tf^{(1)}$ of continuous maps between compact Hausdorff spaces can be naturally identified with the full subcategory of the category $\Tf$ of spaceoids with index set $\Os$ containing a single element.

To every object $X\in \Ob_{\Tf^{(1)}}$ we associate the trivial $\CC$-line bundle $\Xs_X\times\CC$ over the involutive category $\Xs_X:=\Delta_X\times\Rs_{\Os_X}$ with $\Os_X:=\{X\}$ the one point set.

To every continuous map $f:X\to Y$ in $\Tf^{(1)}$ we associate the morphism $(g,\G)$ with $g_\Delta(p,p):=(f(p),f(p))$, $g_\Rs: (X,X)\mapsto (Y,Y)$ and $\G:=\iota_{\Xs_X\times\CC}$.

Note that the trivial bundle over $\Xs_X$ is naturally a $f$-bull-back of the trivial bundle over $\Xs_Y$ and hence $\G$ can be taken as the identity map.

\medskip

Let $\Cs$ and $\Ds$ be two full commutative small \cs-categories (with the same cardinality of the set of objects). Denote by $\Cs_o$ and $\Ds_o$ their sets of identities.

A morphism $\Phi:\Cs\to \Ds$ is an object bijective $*$-functor, i.e.~a map such that
\begin{gather*}
\Phi(x+y)=\Phi(x)+\Phi(y), \ \forall x,y\in \Cs_{AB}, \\
\Phi(a\cdot x)=a\cdot\Phi(x), \ \forall x\in\Cs, \ \forall a\in \CC, \\
\Phi(x\circ y)=\Phi(x)\circ\Phi(y),\  \forall x\in \Cs_{CB},\ y\in \Cs_{BA}\\
\Phi(x^*)=\Phi(x)^*,\ \forall x\in \Cs_{AB}, \\
\Phi(\iota)\in\Ds_o, \ \forall \iota\in\Cs_o, \\
\Phi_o:=\Phi|_{\Cs_o}:\Cs_o\to\Ds_o \quad \text{is bijective}.
\end{gather*}

To every spaceoid $(\Es,\pi,\Xs)$, with $\Xs:=\Delta_X\times\Rs_\O$, we can associate a full commutative \cs-category $\Gamma(\Es)$ as follows:
\begin{itemize}
\item
$\Ob_{\Gamma(\Es)}:=\Os$;
\item
$\forall A,B\in \Ob_{\Gamma(\Es)}$, $\Hom_{\Gamma(\Es)}(B,A):=\Gamma(\Delta_X\times\{(A,B)\};\Es)$,
where $\Gamma(\Delta_X\times\{(A,B)\};\Es)$ denotes the set of continuous sections
$\sigma:\Delta_X\times\{(A,B)\}\to\Es$,
$\sigma:p_{AB}\mapsto \sigma^{AB}_p\in\Es_{p_{AB}}$
of the restriction of $\Es$ to the base space
$\Delta_X\times\{(A,B)\}\subset \Xs$; 
\item
for all $\sigma\in \Hom_{\Gamma(\Es)}(A,B)$ and $\rho\in \Hom_{\Gamma(\Es)}(B,C)$:
\begin{gather*}
\rho\circ\sigma:p_{AC} \mapsto (\rho\circ\sigma)^{AC}_p:=\rho^{AB}_p\circ\sigma^{BC}_p,
\\
\sigma^*: p_{BA} \mapsto (\sigma^*)^{BA}_p:=(\sigma^{AB}_p)^*,
\\
\|\sigma\|:=\sup_{p\in \Delta_X}\|\sigma^{AB}_p\|_{\Es},
\end{gather*}
with operations taken in the total space $\Es$ of the Fell bundle.
\end{itemize}

We extend now the definition of $\Gamma$ to the morphism of $\Tf$ in order to obtain a contravariant functor.

Let $(f,\F)$ be a morphism in $\Tf$ from $(\Es_1,\pi_1,\Xs_1)$ to  $(\Es_2,\pi_2,\Xs_2)$.

Given a section $\sigma\in \Gamma(\Es_2)$, we consider the unique section
$f^\bullet(\sigma):\Xs_1\to f^\bullet(\Es_2)$
such that $f^{\pi_2}\circ f^\bullet(\sigma)=\sigma\circ f$ and
the composition $\F\circ f^\bullet(\sigma)$.

In this way we get a map
\begin{equation*}
\Gamma_{(f,\F)}: \Gamma(\Es_2)\to\Gamma(\Es_1), \quad \Gamma_{(f,\F)}:\sigma\mapsto\F\circ f^\bullet(\sigma), \quad \forall \sigma\in \Gamma(\Es_2).
\end{equation*}

\begin{proposition} (\cite[Proposition~4.1]{BCL17})
Let $(\Es_1,\pi_1,\Xs_1)\xrightarrow{(f,\F)} (\Es_2,\pi_2,\Xs_2)$ be a morphism in $\Tf$, the map $\Gamma_{(f,\F)}:\Gamma(\Es_2)\to \Gamma(\Es_1)$ is a morphism in $\Af$.

The pair of maps $\Gamma:(\Es,\pi,\Xs) \mapsto \Gamma(\Es)$ and $\Gamma:(f,\F)\mapsto\Gamma_{(f,\F)}$ gives a contravariant functor from the category $\Tf$ of spaceoids to the category $\Af$ of small full commutative \cs-categories.
\end{proposition}

We proceed to associate to every commutative full
\cs-category $\Cs$ its spectral spaceoid $\Sigma(\Cs):=(\Es^\Cs,\pi^\Cs,\Xs^\Cs)$, see~\cite[Section~5]{BCL17} for details.

\begin{itemize}
\item
The set $[\Cs;\CC]$ of $\CC$-valued $*$-functors $\omega:\Cs\to\CC$, with the weakest topology making all evaluations continuous, is a compact Hausdorff topological space.
\item
By definition two $*$-functors $\omega_1,\omega_2\in [\Cs;\CC]$ are \emph{unitarily equivalent} if there exists a ``unitary'' natural trasformation $A\mapsto\nu_A\in \TT$ between them.
This is true iff $\omega_1|_{\Cs_{AA}}=\omega_2|_{\Cs_{AA}}$ for all $A\in \Ob_\Cs$.
\item
Let $\Sp_b(\Cs):=\{[\omega] \ | \ \omega\in [\Cs;\CC]\}$ denote the \emph{base spectrum} of $\Cs$, defined as the set of unitary equivalence classes of $*$-functors in $[\Cs;\CC]$.
It is a compact Hausdorff space with the quotient topology induced by the map  $\omega\mapsto[\omega]$.
\end{itemize}

\begin{itemize}
\item
Let $\Xs^\Cs:=\Delta^{\Cs}\times\Rs^{\Cs}$ be the direct product topological involutive category 
of the compact Hausdorff $*$-category $\Delta^\Cs:=\Delta_{\Sp_b(\Cs)}$ and the topologically discrete $*$-category $\Rs^\Cs:=\Cs/\Cs\simeq\Rs_{\Ob_\Cs}$.
\item
For $\omega\in [\Cs;\CC]$, the set $\Is_\omega:=\{x\in \Cs\ | \ \omega(x)=0\}$ is an ideal in $\Cs$ and $\Is_{\omega_1}=\Is_{\omega_2}$ if $[\omega_1]=[\omega_2]$.
\item
Denoting by $[\omega]_{AB}$ the point $([\omega],(A,B))\in \Xs^\Cs$, we define:
\begin{gather*}
\Is_{[\omega]_{AB}}:=\Is_\omega \cap \Cs_{AB}, \quad
\Es^\Cs_{[\omega]_{AB}}:=\frac{\Cs_{AB}}{\Is_{[\omega]_{AB}}}, \quad
\Es^\Cs:=\biguplus_{[\omega]_{AB}\in\Xs^\Cs}\Es^\Cs_{[\omega]_{AB}}.
\end{gather*}
\end{itemize}
\begin{proposition}(\cite[Proposition~5.7]{BCL17})
The map $\pi^\Cs: \Es^\Cs\to \Xs^\Cs$, that sends an element $e\in \Es^\Cs_{[\omega]_{AB}}$ to the point $[\omega]_{AB}\in \Xs^\Cs$
has a natural structure of unital rank-one Fell bundle over the topological involutive category $\Xs^\Cs$.
\end{proposition}

Let $\Phi:\Cs\to\Ds$ be an object-bijective $*$-functor between two small commutative full \cs-categories with spaceoids $\Sigma(\Cs),\Sigma(\Ds)\in \Tf$.

We define a morphism $\Sigma^\Phi:\Sigma(\Ds)\xrightarrow{(\lambda^\Phi,\Lambda^\Phi)}\Sigma(\Cs)$ in the category $\Tf$:
\begin{itemize}
\item
$\lambda^\Phi:\Xs^\Ds\xrightarrow{(\lambda^\Phi_\Delta,\lambda^\Phi_\Rs)}\Xs^\Cs$ where \\
$\lambda^\Phi_\Rs(A,B):=(\Phi_o^{-1}(A),\Phi_o^{-1}(B))$, for all $(A,B)\in \Rs_{\Ob_\Ds}$; \\
$\lambda^\Phi_\Delta([\omega]):=[\omega\circ\Phi] \in \Delta_{\Sp_b(\Cs)}$, for all
$[\omega]\in \Delta_{\Sp_b(\Ds)}$.
\item
The bundle $\biguplus_{[\omega]_{AB}\in\Xs^\Ds} \  \frac{\Cs_{\lambda_\Rs^\Phi(AB)}}{\Is_{\lambda^\Phi([\omega]_{AB})}}$ with the maps \\
$\pi^\Phi:
([\omega]_{AB},\ x+\Is_{\lambda^\Phi([\omega]_{AB})})\mapsto [\omega]_{AB}\in \Xs^\Ds, \quad x\in \Cs_{\lambda^\Phi_\Rs(AB)}$,
\\
$\Phi^\pi: ([\omega]_{AB},\ x+\Is_{\lambda^\Phi([\omega]_{AB})})\mapsto
(\lambda^\Phi([\omega]_{AB}), x+\Is_{\lambda^\Phi([\omega]_{AB})})\in\Es^\Cs
$ \\
is a $\lambda^\Phi$-pull-back $(\lambda^\Phi)^\bullet(\Es^\Cs)$ of the Fell bundle $(\Es^\Cs,\pi^\Cs,\Xs^\Cs)$.
\end{itemize}

\begin{itemize}
\item
Since $\Phi(\Is_{\lambda^\Phi([\omega]_{AB})})\subset \Is_{[\omega]_{AB}}$ for $[\omega]_{AB}\in \Xs^\Ds$, we can define a map  \\
$\Lambda^\Phi: (\lambda^\Phi)^\bullet(\Es^\Cs)
\to \Es^\Ds$ by
$\Big([\omega]_{AB}, \ x+\Is_{\lambda^\Phi([\omega]_{AB})}\Big)
\mapsto \Big([\omega]_{AB}, \ \Phi(x) +\Is_{[\omega]_{AB}}\Big)$.
\end{itemize}

\begin{proposition} (\cite[Proposition~5.8]{BCL17}) 
For any morphism $\Cs\xrightarrow{\Phi}\Ds$ in $\Af$, the map $\Sigma(\Ds)\xrightarrow{\Sigma^\Phi}\Sigma(\Cs)$ is a morphism of spectral spaceoids.
The pair of maps $\Sigma:\Cs\mapsto\Sigma(\Cs)$ and $\Sigma: \Phi\mapsto\Sigma^\Phi$ give a contravariant functor $\Sigma:\Af\to\Tf$, from the category $\Af$ of object-bijective $*$-functors between small commutative full \cs-categories to the category $\Tf$ of spaceoids.
\end{proposition}

We can now state our main duality theorem for commutative full \cs-categories:

\begin{theorem} (P.~Bertozzini-R.~Conti-W.~Lewkeeratiyutkul~\cite[Theorem~6.5]{BCL17})
There exists a duality $(\Gamma,\Sigma)$ between the category $\Tf$ of object-bijective morphisms between topological spaceoids and the category $\Af$ of object-bijective $*$-functors between small commutative full \cs-categories, where
\begin{itemize}
\item
$\Gamma$ is the functor that to every spaceoid $(\Es,\pi,\Xs)\in \Ob_{\Tf}$ associates the small commutative full \cs-category $\Gamma(\Es)$ and that to every morphism between topological spaceoids $(f,\F): (\Es_1,\pi_1,\Xs_1)\to (\Es_2,\pi_2,\Xs_2)$
associates the $*$-functor $\Gamma_{(f,\F)}$;
\item
$\Sigma$ is the functor that to every small commutative full \cs-category $\Cs$ associates its spectral spaceoid $\Sigma(\Cs)$
and that to every object-bijective $*$-functor $\Phi:\Cs\to \Ds$ of \cs-categories in $\Af$ associates
the morphism $\Sigma^\Phi:\Sigma(\Ds)\to\Sigma(\Cs)$
between spaceoids.
\end{itemize}
\end{theorem}
The natural isomorphism $\Gg: \I_\Af \to \Gamma\circ\Sigma$ is provided by the \emph{horizontally categorified Gel'fand transforms} $\Gg_\Cs: \Cs\to\Gamma(\Sigma(\Cs))$
defined by
\begin{gather*}
\Gg_\Cs: \Cs\to\Gamma(\Es^\Cs), \quad
\Gg_\Cs: x\mapsto\hat{x}\quad \text{where} \\
\hat{x}^{AB}_{[\omega]}:=x+\Is_{[\omega]_{AB}}, \quad \forall x\in \Cs_{AB}.
\end{gather*}
\begin{proposition} (\cite[Theorem~6.3]{BCL17})
The functor $\Gamma: \Tf\to \Af$ is representative i.e.~given a commutative full \cs-category $\Cs$, the Gel'fand transform
$\Gg_\Cs :\Cs\to\Gamma(\Sigma(\Cs))$ is a full isometric (hence faithful) $*$-functor.
\end{proposition}

The natural isomorphism $\Eg: \I_\Tf \to \Sigma\circ\Gamma$ is provided by the \emph{horizontally categorified ``evaluation'' transforms}
$\Eg_\Es: (\Es,\pi,\Xs)\xrightarrow{(\eta^\Es,\Omega^\Es)}\Sigma(\Gamma(\Es))$,
defined as follows:
\begin{itemize}
\item
$\eta^\Es_\Rs (A,B):=(A,B), \quad \forall (A,B)\in \Rs_\Os$.
\item
$\eta^\Es_\Delta: \Delta_X\to \Delta_{\Sp_b(\Gamma(\Es))}$, 
$p\mapsto [\gamma\circ\ev_p]$, where $\ev_p: \Gamma(\Es)\to \uplus_{(AB)\in\Rs_\O} \ \Es_{p_{AB}}$ is the evaluation map
given by
$\sigma\mapsto \sigma^{AB}_p$ that is a $*$-functor with values in a one dimensional \cs-category that actually determines\footnote{There is always a $\CC$ valued $*$-functor $\gamma:\uplus_{(AB)\in\Rs_\O} \ \Es_{p_{AB}}\to \CC$ and any two compositions of $\ev_p$ with such $*$-functors are unitarily equivalent because they coincide on the diagonal \cs-algebras $\Es_{p_{AA}}$.} a unique point
$[\gamma\circ\ev_p]\in \Delta_{\Sp_b(\Gamma(\Es))}$.
\item
$\biguplus_{p_{AB}\in\Xs} \  \Gamma(\Es)_{\eta^\Es_\Rs(AB)}/\Is_{\eta^\Es(p_{AB})}$ 
with the projection
$(p_{AB}, \ \sigma + \Is_{\eta^\Es(p_{AB})})\mapsto p_{AB}$, and with the $\Es^{\Gamma(\Es)}$-valued map
$(p_{AB}, \ \sigma + \Is_{\eta^\Es(p_{AB})}) \mapsto \sigma + \Is_{\eta^\Es(p_{AB})}$,
is a $\eta^\Es$-pull-back $(\eta^\Es)^\bullet(\Es^{\Gamma(\Es)})$ of $\Sigma(\Gamma(\Es))$.
\item
$\Omega^\Es:(\eta^\Es)^\bullet(\Es^{\Gamma(\Es)})\to \Es$ is defined by  \\
$\Omega^\Es: (p_{AB}, \ \sigma+\Is_{\eta^\Es(p_{AB})})
\mapsto \sigma^{AB}_p,
\quad \forall \sigma \in \Gamma(\Es)_{AB}, \quad p_{AB}\in\Xs$.
\end{itemize}

In particular, with such definitions we can prove:
\begin{proposition} (\cite[Theorem~6.4]{BCL17})
The functor $\Sigma: \Af \to \Tf$ is representative i.e.~given a spaceoid $(\Es,\pi,\Xs)$, the evaluation transform
$\Eg_\Es :(\Es,\pi,\Xs)\to\Sigma(\Gamma(\Es))$ is an isomorphism in the category of spaceoids.
\end{proposition}

\medskip

We are now working on a number of generalizations and extensions of our horizontal categorified Gel'fand duality:
\begin{itemize}
\item[$\looparrowright$]
The first immediate possibility is to extend Gel'fand duality to include the case of categories of general $*$-functors between full commutative \cs-categories. This will necessarily require the consideration of categories of $*$-relators (see~\cite{BCL6}) between \cs-categories.
\item[$\looparrowright$]
Our duality theorem is for now limited to the case of full commutative \cs-categories and further work is necessary in order to extend the result to a Gel'fand duality for non-full \cs-categories.
\item[$\looparrowright$] 
Very interesting is the possibility to generalize our duality to a full spectral theory for non-commutative \cs-categories and Fell bundles in term of endofunctors with target in the category of Fell line-bundles. This might be a step in order to make contact with the notion  of ``Fell bundle geometry'' introduced by R.~Martins~\cite{Marti1,Marti2,Marti3} for a categorical reformulation of the spectral triple in the standard model. 
\item[$\looparrowright$] 
Furthermore we would like to explore if our approach will allow to develop categorifications of J.~Dauns-K.H.~Hofmann 
theorem~\cite{DH}, R.~Cirelli-A.~Mani\`a-L.~Pizzocchero~\cite{CMP} spectral theorem and 
G.~Elliott-K.~Kawamura~\cite{EK, Kaw} Serre-Swan equivalence for general non-commutative C*-algebras. 
\item[$\looparrowright$] 
Similarly, it might be important to study the relation between our spectral spaceoids and other spectral notions such
as locales and topoi that are already used in the constructive spectral theorems by B.~Banachewki-C.~Mulvey~\cite{BM} and 
C.~Heunen-K.~Landsman-B.~Spitters~\cite{HS,HLS3,HLS4}. 
In the same order of ideas, motivated by a general spectral theory for \cs-categories, it is worth investigating in the non-commutative case the connection between \cs-categories,  spectral spaceoids and categorified notions of (locale) quantale already developed for (commutative) \cs-algebras (see D.~Kruml-J.~Pelletier-P.~Resende-J.~Rosicky~\cite{KPRR}, D.~Kruml-P.~Resende~\cite{KrR} P.~Resende~\cite{Res}, L.~Crane~\cite{Cr2} and references therein for details).
\item[$\looparrowright$]
The existence of a horizontal categorified Gel'fand transform might be relevant for the study of harmonic analysis on commutative groupoids. In this direction it is natural to investigate the implications for a Pontrjagin duality for commutative groupoids 
and later, in a fully non-commutative context, the relations with the theory of \cs-pseudo-multiplicative unitaries that has been recently developed by T.~Timmermann~\cite{Ti,Ti2,Ti3,Ti4}.
\item[$\looparrowright$]
Extremely intriguing for its possible physical implications in algebraic quantum field theory is the appearance of a natural ``local gauge structure'' on the spectra: the spectrum is no more just a (topological) space, but a special fiber bundle.
Possible relations with the work of E.~Vasselli~\cite{Va0,Va,Va2,Va3} on continous fields of \hbox{\cs-categories} in the theory of superselection sectors and especially with the recent work on net bundles and gauge theory by J.~Roberts-G.~Ruzzi-E.~Vasselli~\cite{RRV,RRV2} remain to be explored.
\end{itemize}

			\subsubsection{Higher \cs-categories.}

In our last forthcoming work\footnote{P.~Bertozzini, R.~Conti, W.~Lewkeeratiyutkul,  N.~Suthichitranont, Strict Higher 
\cs-categories, in preparation.}, we proceed to further extend the categorification process of Gel'fand duality theorem to a full ``vertical categorification''~\cite{Ba1}.

For this purpose we first provide, via globular sets (see T.~Leinster's book~\cite{Le}), a suitable definition of ``strict'' $n$-\cs-category.

In practice, without entering here in further technical details (see the slides~\cite[Pages~93-104]{B2} for a deeper overview) a strict higher \cs-category $\Cs$ (or more generally a higher Fell bundle over a higher inverse $*$-category $\Xs$), is provided by a strict higher $*$-category $\Cs$ fibered over a strict higher inverse 
$*$-category $\Xs$ whose compositions and involutions satisfy, fiberwise at all levels, ``appropriate versions'' of all the properties listed in the definition of a Fell bundle.

In the special case of commutative full strict $n$-\cs-categories, we develope a spectral Gel'fand theorem in term of $n$-spaceoids i.e.~rank-one $n$-\cs-Fell bundles over a ``particular'' $n$-$*$-category (that is given by the direct product of the diagonal equivalence relation of a compact Hausdorff space and the quotient $n$-$*$-category $\Cs/\Cs$ of an $n$-\cs-category $\Cs$).

\begin{itemize}
\item[$\looparrowright$]
Unfortunately our definition is for now limited to the case of strict higher \hbox{\cs-cat}\-egories. Of course, as always the case in higher category theory, an even more  interesting problem will be the characterization of suitable axioms for ``weak higher 
\cs-categories''.\footnote{For the purpose of the development of a notion of ``weak C*-algebra'' (where the usual axioms for product, identity and involution are expected to hold only up to isomorphism) it is interesting to consider the recent work by 
P.~Bouwknegt-K.~Hannabuss-V.~Mathai~\cite{BHM}, where ``C*-algebras'' with a strictly non-associative product are defined as objects internal to suitable monoidal dagger categories.} 
This is one of the main obstacles in the development of a full categorification of the notion of spectral triple and of A.~Connes non-commutative geometry.
\item[$\looparrowright$]
Note that several examples and definitions of 2-\cs-categories are already available in the literature (see for example 
R.~Longo-J.~Roberts~\cite{LR} and P.~Zito~\cite{Z}).
In general such cases will not exactly fit with the strict version of our axioms for \hbox{$n$-\cs-cat}egories. Actually we expect to have a complete hierarchy of definitions of higher \cs-categories depending on the ``depth'' at which some axioms are required to be satisfied (i.e.~some properties can be required to hold only for $p$-arrows with $p$ higher than a certain depth).
\item[$\looparrowright$]
In our work, we define (Hilbert \cs-)modules over strict $n$-\cs-categories and in this way we can provide interesting definitions of 
$n$-Hilbert spaces and start a development of ``higher functional analysis''. Extremely interesting for us will be to understand the relation with categorified notions of higher vector and Hilbert spaces developed by M.~Kapranov- V.~Voevodsky~\cite{KV}, 
J.~Baez-A.~Crans~\cite{Ba0,BC}, J.~Elgueta~\cite{E} and J.~Morton~\cite{Mor3,Mor4}. 
\end{itemize}

	\subsection[Categorical NCG and Topoi]{Categorical Non-commutative Geometry and \\ Non-commutative Topoi.}

One of the main goals of our investigation is to discuss the interplay between ideas of categorification and non-commutative geometry. Here there is still much to be done and we can present only a few suggestions. Work is in progress.
\begin{itemize}
\item[$\looparrowright$]
Every isomorphism class of a full commutative \cs-category can be identified with an equivalence relation in the Picard-Morita $1$-category of Abelian unital \cs-algebras.
In practice a \cs-category is just a ``strict implementation'' of an equivalence relation subcategory of Picard-Morita.

Since morphism of spectral triples (more generally morphisms of non-commutative spaces)
are essentially ``special cases'' of Morita morphisms, we started the study of ``spectral triples over \cs-categories''
and
we are now trying to develop a notion of horizontal categorification of spectral triples (and of other spectral geometries) in order to identify a correct definition of morphism of spectral triples that supports a duality with a suitable spectrum (in the commutative case).

The general picture that is emerging~\cite[Pages~105-108]{B2}\footnote{P.~Bertozzini-R.~Conti-W.~Lewkeeratiyutkul, Spectral Geometries over C*-categories and Morphisms of Spectral Geometries, in preparation; 
Horizontal Categorification of Spectral Triples, work in progress.} is that a correct notion of metric morphism between spectral triples is given by a kind of ``bivariant version'' of spectral triple i.e.~a bimodule over two different algebras that is equipped with a  left/right action of ``Dirac-like'' operators.
\item[$\looparrowright$]
As a very first step in the direction of a full ``higher non-commutative geometry''\footnote{On this topic the reader is strongly advised to read the interesting discussions on the ``$n$-category caf\'e''~\hlink{http://golem.ph.utexas.edu/category/}{http://golem.ph.utexas.edu/category/} and in particular: U.~Schreiber, Connes Spectral Geometry and the Standard Model II, 06 September 2006.} we plan to start the study of a strict version of ``higher spectral triples'' i.e.~spectral triples over strict higher \cs-categories. As in the case of horizontal categorification, this will provide some hints for a correct definition of ``higher spectral triples''.
\item[$\looparrowright$]
Although at the moment it is only a speculative idea, it is very interesting to explore the possible relation between such ``higher spectra'' (higher spaceoids) and the notions of stacks and gerbes already used in higher gauge theory. The recent work by C.~Daenzer~\cite{Dae} in the context of T-duality discuss a Pontryagin duality between commutative principal bundles and gerbes that might be connected with our categorified Gel'fand transform for commutative \cs-categories.
\item[$\looparrowright$]
Extremely intriguing is the possible connection between the notions of (category of) spectral triples and A.~Grothendieck topoi.
Speculations in this direction have been given by P.~Cartier~\cite{Car} and are also discussed by A.~Connes~\cite{C14}.
A full (categorical) notion of non-commutative space (non-commutative Klein program / non-commutative Grothendieck topos) is still waiting to be defined\footnote{P.~Bertozzini, R.~Conti, W.~Lewkeeratiyutkul, Non-commutative Klein-Cartan Program, work in progress.}.
\end{itemize}

Actually some interesting proposal for a definition of a ``quantum topos'' is already available in the recent work by L.~Crane~\cite{Cr2} based on the notion of ``quantaloids'', a categorification of the notion of quantale (see P.~Resende~\cite{Res} and references therein).

At this level of generality, it is important to emphasize that our discussion of non-commutative geometry has been essentially confined to the consideration of A.~Connes'  approach.
In the field of algebraic geometry (see V.~Ginzburg~\cite{Gi}, M.~Kontsevich-Y.~Soibelman~\cite{KS,KS2} and S.~Mahanta~\cite{Mah1,Mah2,Mah3,Mah4} as recent references), many other people have been trying to propose definitions of non-commutative schemes and non-commutative spaces (see for example A.~Rosenberg~\cite{R} and M.~Kontsevich-A.~Rosenberg~\cite{KR}) as ``spectra'' of Abelian categories (or generalization of Abelian categories such as triangulated, dg, or $A_\infty$ categories). Since every Abelian category is essentially a category of modules, it is in fact usually assumed that an Abelian category should be considered as a topos of sheaves over a non-commutative space.
\begin{itemize}
\item[$\looparrowright$]
It is worth noting that the categories naturally arising in the theory of self-adjoint operator algebras and in A.~Connes' non-commutative geometry are $*$-monoidal categories (see~\cite{BCL17} for detailed definitions). The monoidal property is perfectly in line with the recent proposal by T.~Maszczyk~\cite{Mas} to construct a theory of algebraic non-commutative geometry based on Abelian categories equipped with a monoidal structure.

At this point it is actually tempting (in our opinion) to think that also the involutive structures (and other properties strictly related to the existence of an involution  including modular theory\footnote{See section~\ref{sec: modular} for some references.}) are going to play some vital role in the correct definition of a non-commutative generalization of space. But this is still speculation in progress!

\item[$\looparrowright$]
Finally, there are strong indications (V.~Dolgushev-D.~Tamarkin-B.~Tsygan~\cite{DTT})\footnote{See also the very detailed discussion on the blog ``$n$-category caf\'e'': J.~Baez, Infinitely Categorified Calculus, 09 February 2007.} coming again from ``algebraic non-commutative geometry'' that a proper categorification of non-commutative geometry might actually be possible only considering $\infty$-categories. The implications for a program of categorification of A.~Connes' spectral triples is not yet clear to us.
\end{itemize}

\section{Applications to Physics.}\label{sec: physics}

In this final section we would like to spend some time to introduce (in a non-technical way) the mathematical readers to the consideration of some extremely important topics in quantum physics that are essentially motivating the construction of non-commutative spaces, the use of categorical ideas and the eventual merging of these two lines of thought.

The two main subjects of our discussion, non-commutative geometry and category theory, have been separately used and applied in theoretical physics (although not as widely as we would have liked to see) and we are going to review here some of the main historical steps in these directions.

Anyway, our feeling is that the most important input to physics will come from a kind of  ``combined'' approach where non-commutative and categorical structures are applied in a ``synergic way'' in an ``algebraic theory of quantum gravity'' (AQG).
A concrete proposal in this direction is presented in section~\ref{sec: AQG}.

	\subsection{Categories in Physics.}

Category theory has been conceived as a tool to formalize basic structures (functors, natural transformations) that are omnipresent in algebraic topology. Its level of abstraction has been an obstacle to its utilization even in the mathematics community and so it does not come as a surprise that fruitful applications to physics had to wait.

Probably, the first person to call for the usage of categorical methods in physics has been  J.~Roberts in the seventies.
The joint work of S.~Doplicher and J.~Roberts~\cite{DR1,DR2} on the theory of superselection sectors in algebraic quantum field theory\footnote{The texts by R.~Haag~\cite{H},  H.~Araki~\cite{A}, D.~Kastler~\cite{K0} and the recent book~\cite{BBIM} contain detailed introductions to superselection theory in algebraic quantum field theory.} is one of the most eloquent examples of the power of category theory when applied to fundamental physics: giving a full explanation of the origin of compact gauge groups of the first kind and field algebras in quantum field theory and providing at the same time a general Tannaka-Kre\u\i n duality theory for compact groups, where the dual of a compact group is given by a particular monoidal W*-category.
Since then, monoidal $*$-categories are a common topic of investigation in algebraic quantum field theory,\footnote{For a complete list of all relevant papers and a recent ``philosophical'' overview of the subject see H.~Halvorson-M.~M\"uger~\cite{HM} and also R.~Brunetti-K.~Fredenhagen~\cite{BrF}.} where several people are still working on possible variants and extensions of superselection theory.\footnote{A large literature is of course available on monoidal categories, tensor categories (see for example M.~M\"uger~\cite{Mu} for a survey) and their application to the theory of ``quantum groups'' (see for example R.~Street~\cite{St} for a clear introduction) as well as many other different subjects, but it is outside the scope of this survey to enter into further details on these topics. 
}

\medskip

The role of categories in physics, more recently, has been stressed also from different areas of research such as conformal field theory (G.~Segal~\cite{Se}) and topological quantum field theory (M.~Atiyah~\cite{At}).
A very interesting relation between axiomatizations of these topological quantum field theories (in their ``extended'' functorial version~\cite{BD1}) and the Haag-Kastler axioms for algebraic quantum field theory has been proposed in the recent work by U.~Schreiber~\cite{Sch}.  

\medskip

C.~Isham has been the pioneer in suggesting to consider topoi as basic structures for the construction of alternative quantum theories in which ordinary set theoretic concepts (including real/complex numbers and classical two valued logic) are replaced by more general topos theoretic notions. His research with J.~Butterfield~\cite{BI1,BI2,BI3,BI4,BI5,BI6,BHI} and more recently with 
A.~D\"oring~\cite{DI1,DI2,DI3,DI4,DI5,Dor1,Dor2} has polarized the attention towards a possible usage of topos theory in quantum mechanics~\footnote{See the recent papers by C.~Heunen-K.~Landsman-B.~Spitters~\cite{HS,HLS,HLS3,HLS4} and the preprint by C.~Flori~\cite{Fl}.} and quantum gravity, an idea that has influenced several other authors working on quantum gravity.

\medskip

S.~Abramsky-B.~Coecke~\cite{AbC0,AbC,AbC2,AbC3,Ab0,Ab3,Co0,Co,Co2,Co3,Co5,Ab4} with their collaborators R.~Duncan-B.~Edwards-D.~Pavlovic-E.~Paquette-S.~Perdrix-J.~Vicary~\cite{AbD,CPav,CPa,CPa2,CPP,CPP2,CD,CD2,CPPe,CE,CPV,Vi1,Vi2} are actively developing a categorical axiomatization for quantum mechanics based on symmetric monoidal categories with intriguing links to knot theory, logic  and computer science~\cite{Ab1,Ab2,Ab6,Ab7}. 

Related works on categorical quantum logic in the setting of dagger-monoidal or dagger categories with kernels are carried on by C.~Heunen-B.~Jacobs~\cite{He1,He2,HeJ}. 

\medskip

N.~Landsman~\cite{La,La2,La3} in his study of quantization and of the relation between classical Poisson geometry and operator algebras of quantum systems, has been constantly exploiting techniques from category theory (groupoids, Morita equivalence).
His recent works M.~Caspers-C.~Heunen-N.~Landsman-B.~Spitters~\cite{HS,CHLS,HLS3,HLS4}, further elaborate on the C.~Isham-J.~Butterfield-A.~D\"oring proposal to base physics on topos theory, opening the way to reconsider algebraic quantum theory as a ``classical theory'' living in a suitable ``spectral topos'' and proposing an extension of general covariance in terms of geometric morphisms between topoi~\cite{HLS}. 

\medskip

J.~Baez has been one of the first pioneers and the most prominent advocate in the development, with J.~Dolan, 
of higher categorical structures~\cite{BD1}  (``opetopic'' $n$-categories~\cite{BD2,Ba1}, categorification~\cite{BD,BD4}, 
$2$-Hilbert spaces~\cite{Ba0}) and in the usage of categorical methods in quantum physics and in quantum 
gravity~\cite{Ba2,Ba3,Ba4}.
J.Baez and his collaborators and students, 
T.~Bartles-A.~Crans-A.~Hoffnung- L.~Langford-A.~Lauda-J.~Morton-M.~Neuchl-C.~Rogers-U.~Schreiber-M.~Shulman-M.~Stay-D.~Stevenson-C.~Walker
have been eleborating huge portions of ``higher algebra extensions'' of mathematics~\cite{BSh}  
(braided monoidal $2$-categories~\cite{BN}, $2$-Lie algebras~\cite{BC}, $2$-tangles~\cite{BLan1,BLan2}, 
$2$-groups~\cite{BLa,BSt,BBFW}, $2$-bundles~\cite{Bart}, groupoidification~\cite{BHW}) 
and their applications to physics: higher gauge theories~\cite{Ba5,Bart,BSc1,BSc2,BCSS}, 
higher symplectic geometry~\cite{BHR,BRo}, quantum computation~\cite{BS} and  ``combinatorial'' quantum 
mechanics~\cite{Mor1}. 




\medskip

A new emerging field of ``categorical quantum gravity'' is developing (see the works by 
L.~Crane~\cite{Cr,Cr2,Cr3,Cr4}, J.~Baez~\cite{Ba2,Ba3,Ba4} and, for a categorical  approach via topological quantum field theory  cobordism, also J.~Morton~\cite{Mor2,Mor3,Mor4}.

	\subsection{Categorical Covariance.}

Covariance of physical theories has been always discussed in the limited domain of groups acting on spaces.
\begin{itemize}
\item
Aristotles' physics is based on the covariance group $SO_3(\RR)$ of rotations in $\RR^3$ that was the supposed symmetry group of a three dimensional vector space with the center of the Earth at the origin.
\item
Galilei's relativity principle requires as covariance group the Galilei group, which is the ten parameters symmetry group of the Newtonian space-time (i.e.~a family $E_t$ of three dimensional Euclidean spaces parametrized by elements $t$ in a one dimensional Euclidean space $T$) generated by 3 space translations, 1 time translation, 3 rotations and 3 boosts.
\item
Poincar\'e covariance group consists of the semidirect product of Lorentz group $\L$ with the group of translations in $\RR^4$ and it is the symmetry group of the four dimensional Minkowski space (an affine four dimensional space modeled on $\RR^4$ with metric of signature $(-+++)$).
\item
Einstein covariance group is the group of diffeomorphisms of a four dimensional Lorentzian manifold (note that in this case the metric and the causal structure is not preserved).
\end{itemize}
Different observers are ``related'' through transformations in the given covariance group.

\medskip

\begin{itemize}
\item[$\looparrowright$]
There is no deep physical or operational reason to think that only groups (or quantum groups) might be the right mathematical structure to capture the ``translation'' between different observers and actually, in our opinion, categories provide a much more suitable environment in which also the discussion of ``partial translations'' between observers can be described. Work is in progress on these issues~\cite{B}.
\end{itemize}

The substitution of groups with categories (or graphs), as the basic covariance structure of theories,  should be a key ingredient for all the approaches based on deduction of physics from operationally founded principles of information theory (see C.~Rovelli~\cite{Ro3} and A.~Grinbaum~\cite{Gri1,Gri2,Gri3}) and, in the context of quantum gravity, also for theories based on the formalism of quantum casual histories (see for example E.~Hawkins-H.~Sahlmann-F.Makopoulou~\cite{HMS} and F.~Markopoulou~\cite{Mar3}).

\medskip

As an example of the relevance of the idea of categorical covariance, we mention several new works by R.~Brunetti-K.~Fredenhagen-R.~Verch~\cite{BFV}, R.~Brunetti-G.~Ruzzi~\cite{BrR} and R.~Brunetti-M.~Porrmann-G.~Ruzzi~\cite{BPR} that,
following the fundamental idea of J.~Dimock~\cite{Dim1, Dim2},
aim at a generalization of  H.~Araki-R.~Haag-D.~Kastler algebraic quantum field theory axiomatization\footnote{See H.~Araki's and R.~Haag's books~\cite{A,H} and also  K.~Fredenhagen-K.-H.~Rehren-H.~Seiler~\cite{FRS} for a discussion and contextualization of algebraic quantum field theory.}, that is suitable for an Einstein covariant background. Similar ideas are also used in the non-commutative versions of the axioms recently proposed by M.~Paschke and R.~Verch~\cite{PV,PV2}.

	\subsection{Non-commutative Space-Time.}\label{sec: nc-st}

There are three main reasons for the introduction of non-commutative space-time structures in physics and for the deep interest developed by physicists for ``non-commutative geometry'' (not only A.~Connes'one):
\begin{itemize}
\item
The awareness that quantum effects (Heisenberg uncertainty principle), coupled to the general relativistic effect of the energy-momentum tensor on the curvature of space-time (Einstein equation), entail that at very small scales the space-time manifold structure might be ``unphysical''.
\item
The belief that modification to the short scale structure of space-time might help to resolve the problems of ``ultraviolet divergences'' in quantum field theory (that arise, by Heisenberg uncertainty, from the arbitrary high momentum associated with arbitrary small length scales) and of ``singularities'' in general relativity.
\item
The intuition that in order to include the remaining physical forces
(nuclear and electromagnetic) in a ``geometrization'' program, going beyond the one realized for gravity by A.~Einstein's general relativity, it might be necessary to make use of geometrical environments more sophisticated than those provided by usual Riemannian/Lorentzian geometry.
\end{itemize}

The first one to conjecture that, at small scales, space-time modeled by ``manifolds'' might not be an operationally defined concept was B.~Riemann himself.
A.~Einstein immediately recognized the need to introduce ``quantum'' modifications to general relativity
and M.~Bronstein realized that the specific problems posed by a covariant quantization of general relativity were calling for a rejection of the usual space-time modeled via Riemannian geometry.
Recently a more complete argument has been put forward by S.~Doplicher-K.~Fredenhagen-J.~Roberts~\cite{DFR0,DFR} and by many other in several variants.

J.~Wheeler~\cite{Wh3} introduced the well-known ``space-time-foam'' term to define the hypothetical geometrical structure that should supersede smooth differentiable manifolds at small scales.
Non-commutative geometries are a natural candidate to replace ordinary Lorentzian smooth manifolds as the arena of physics and provide a rigorous (although incomplete, yet) formalization of the notion of space-time ``fuzziness''.

\medskip

The notion of non-commutative space-time originated from an idea of  W.~Heisenberg~\footnote{He communicated it in a letter to R.~Peierls who shared the suggestion with W.~Pauli and R.~Oppenheimer.} that was developed by H.~Snyder~\cite{Sn}.
More recently  S.~Doplicher-K.~Fredenhagen-J.~Roberts~\cite{DFR0,DFR,Do1,Do2} described a new version of Poincar\'e covariant non-commutative space.
An algebraic quantum field theory on such non-commutative spaces is currently under active development (see S.~Doplicher~\cite{Do} for a recent review) and there are some hopes to get in such cases a theory that is free from divergences.
Many other ``variants'' of non-commutative space-time (mostly obtained by ``deformation'' of Minkowski space-time or as ``homogeneuos spaces'' of a ``deformed'' Poincar\'e group) and non-commutative field theory on them are now under investigation in theoretical physics (see for example J.~Madore~\cite{Ma}, B.~Cerchiai-G.~Fiore-J.~Madore~\cite{CFM}, G.~Fiore~\cite{Fio,Fio2,Fio3}, G.~Fiore-J.Wess~\cite{FW}, H.~Grosse-G.~Lechner~\cite{GLe} and references therein), but it is beyond our scope here to enter the details of their description.

\medskip

On the path of complete ``geometrization of matter'' envisioned by B.~Riemann and W.~K.~Clifford, A.~Einstein has been one of the few to stress the conceptual need for a geometrical treatment of the nuclear and electromagnetic forces alongside with gravity.
T.~Kaluza and O.~Klein's theory of unification of electromagnetism with gravity via ``extra-dimensional'' Lorentzian manifolds was clearly going in this direction, but it has gained some popularity only recently, with the introduction of superstrings that, for reasons of internal consistency, require the existence of (compactified) extra dimensions and whose treatment of gravity is manifestly non-background-independent (in the sense required by general relativity).

To date, the most successful achievement in the direction of ``geometrization of physical interactions'', has been obtained by A.~Connes~\cite{C7,CLo,C,CM6} (see also the works by  A.~Chamseddine-A.~Connes~\cite{CC,CC2,CC3,CC4,Ch},  
A.~Chamseddine-A.~Connes-M.~Marcolli~\cite{CCMa} and J.~Barrett~\cite{Bar} for a Lorentzian version) who has promoted the view that the complexity of the standard model in particle physics should be reconsidered as revealing the features of the non-commutative geometry of space-time. The program describes (for the moment only at the classical level) how gravity and all the other fundamental interactions of particle physics arise as a kind of gravitational field on a non-commutative space-time given by a spectral triple over a \cs-algebra that is a tensor product of the algebra of continous function on a $4$-dimensional orientable spin-manifold and a finite dimensional real \cs-algebra.

	\subsection{Spectral Space-Time.}

What we call here ``spectral space-time'' is the idea that space-time (commutative or not) has to be ``reconstructed a posteriori'', from other operationally defined degrees of freedom, in a spectral way. The origin of such ``pregeometrical philosophy'' is less clear.

Space-time as a ``relational'' a posteriori entity originate from ideas of G.W.~Leibniz, G.~Berkeley, E.~Mach.

Although pregeometrical speculations, in western philosophy, probably date as far back as Pythagoras, their first modern incarnation probably starts with J.~Wheeler's ``pregeometries''~\cite{Wh1,MLP} and ``it from bit''~\cite{Wh2} proposals.

R.~Geroch~\cite{Ge}, with his Einstein algebras,
was the first to suggest a ``transition'' from spaces to algebras in order to solve the problem of ``singularities'' in general relativity.


The fundamental idea that space-time can be recovered from the specification of suitable states of the system, has been the subject of scattered speculations in algebraic quantum field theory in the past by A.~Ocneanu~\footnote{As reported in A.~Jadczyk~\cite{Ja}.}, S.~Doplicher~\cite{Do0}, U.~Bannier~\cite{Ban} and, in the ``modular localization program'' (see R.~Brunetti-D.~Guido-R.~Longo~\cite{BGL} and references therein), has been conjectured by 
N.~Pinamonti~\cite{Pi}.

Extremely important rigorous results including a complete reconstruction of Minkowski space-time~\cite{SuW} have been achieved in the ``geometric modular action'' program by D.~Buchholz-S.~J.~Summers (see D.~Buchholz-S.~J.~Summers~\cite{BuS,BS2}, D.~Buchholz-M.~Florig-S.~J.~Summers~\cite{BFS},  D.~Buchholz-O.~Dreyer-M.~Florig-S.~J.~Summers~\cite{BDFS}, for details and S.~J.~Summers~\cite{Su2} for an excellent review and additional references).

More recently the idea has gained importance in the light of attempts to reconstruct quantum physics from operationally founded  quantum information (among others, J.~Bub-R.~Clifton-H.~Halvorson~\cite{BCH}, A.~Grinbaum~\cite{Gri1,Gri2,Gri3,Gri4} and especially C.~Rovelli's suggestion~\cite[section~5.6.4]{Ro}), but in its full generality, the recostruction of space-time is still an unsolved problem.

\begin{itemize}
\item[$\looparrowright$]
This is probably because only now the Araki-Haag-Kastler axiomatization
has been suitably extended to incorporate general covariance  (R.~Brunetti-K.~Fredenhagen-R.~Verch~\cite{BFV}), but there are, in our opinion, other fundamental issues that need to be addressed in a completely unconventional way and that are related to the ``philosophical interpretation'' of states and observables in the theory in ``atemporal-covariant'' context (following ideas of C.~Isham and collaborators~\cite{I3,I4,IL0,IL1,ILSS}, C.~Rovelli and collaborators~\cite{Ro,Ro3,RR,MPR}, J.~Hartle~\cite{Hart}, L.~Hardy~\cite{Har1,Har2,Har3}, J.~Dowling-S.~Jay Olson~\cite{DJO1,DJO2}).
\end{itemize}

That essential information about the underlying space-time is already contained in the algebra of observables of the system (and its Hilbert space representation) is clearly indicated by R.~Feynman-F.~Dyson~\cite{Dy} reconstruction of Maxwell equations (and hence of the Poincar\'e group of symmetries) from the commutation relations of ordinary non-relativistic quantum mechanics of a free particle, an argument recently revised
and extended to non-commutative configuration spaces by T.~Kopf-M.~Paschke~\cite{P0,KP3}.

\medskip

In a slightly different context, in their discussion of the construction of the quantum theory of spin particles on a (compact Riemannian manifold),  J.~Fr\"ohlich-O.~Grandjean-A.~Recknagel~\cite{FGR,FGR2,FGR3,FGR4}, have been considering  several important unsolved aspects of the relationship between the underlying configuration space of a physical system and the actual non-commutative geometry exhibited at the level of its algebra of observables (phase-space). The solution of these problems is still fundamental in the construction of a theory of spectral space-time and quantum gravity based on algebras of observables and their states. We will have more to say about this problem in the final section~\ref{sec: AQG}.

\medskip

\begin{itemize}
\item[$\looparrowright$]
That non-commutative geometry provides a suitable environment for the implementation of spectral reconstruction of space-time from states and observables in quantum physics has been the main motivating idea of one us (P.B.) since 1990 and it is still an open work in progress~\cite{B1}.
\end{itemize}

	\subsection{Quantum Gravity.} \label{sec: qg}

Quantum gravity is the discipline of theoretical physics that deals with the interplay between quantum physics and general relativity. The need for research in this direction was actually recognized by A.~Einstein since the birth of general relativity and several people started to work on it from 1930. Unfortunately, after many years of research by some of the best scientists, we do not have yet an established theory, let alone a mathematically sound frame for these questions.

Following closely C.~Isham's excellent reviews~\cite{I2,I5}, here below we try to summarize the several  approaches to quantum gravity:\footnote{See also Appendix~C in C.~Rovelli's book~\cite{Ro} for a detailed history of the subject.}
\begin{itemize}
\item[a)] \emph{Quantizations of general relativity.} \\
Approaches of this kind, try to make use of a ``standard version'' of quantum mechanics to substitute (a modified) general relativity with a quantized version.
\begin{itemize}
\item
\emph{Canonical quantization} (initiated by P.~Dirac-P.~Bergmann,
developed by R.~Arnowitt-S.~Deser-C.~Misner and J.~Wheeler-B.~DeWitt and recently revived by A.~Sen-A.~Ashtekar and L.~Smolin-L.~Crane-C.~Rovelli-R.~Gambini and others) is probably the first non-perturbative proposal: it tries to find suitable canonical variables to describe the dynamics of classical general relativity and to perform a quantization on them.
After a period of stagnation, this approach has been revived under the name of
\emph{loop quantum gravity} and it is currently the most elaborate non-perturbative (background-independent) program in quantum gravity (see C.~Rovelli~\cite{Ro,Ro4} for an introduction and also T.~Thiemann~\cite{Th1,Th2,Th3,Th4}).
\item
\emph{Covariant quantization} (initiated by
L.~Rosenfeld,  M.~Fierz-W.~Pauli and developed
by B.~DeWitt, R.~Feynman, G.~'t Hooft)
is a background-dependent perturbative approach (usually the one preferred by particle physicists) in which the non-Minkowskian part of the metric tensor is considered as a classical field propagating on a fixed Minkowski space and quantized as any other such field.
The proof of non-renormalizability of general relativity in this setting has somehow stopped any further attempts in this direction forcing researchers to take the stand that general relativity is not a fundamental theory and prompting the development of supergravity and later string theory (see the approaches listed in item~c) here below). 
This approach is also receiving renewed attention because of important results in the \emph{asymptotic safety scenario} originally suggested by S.~Weinberg~\cite{Wei} and developed by M.~Niedermaier-M.~Reuter~\cite{NR} and R.~Percacci~\cite{Per}. 
\item
\emph{Path integral quantization} (initiated by C.~Misner-J.~Wheeler, developed by S.~Hawking-J.~Hartle)
is a non-perturbative proposal that is characterized by its use of the formalism of
Feynman functional integrals for quantization. In its first incarnation, \emph{Euclidean quantum gravity}, the theory was performing a path quantization of a Riemannian version of general relativity and it was motivated by semiclassical studies by S.~Hawking on the thermodynamic properties of black-holes (quantum field theory on curved space-times). 

Discretized versions of functional integral quantization (see~\cite{Lol} for a review) have been originally based on \emph{Regge calculus} proposed by T.~Regge~\cite{Reg}, but  recently the approach has been revived in a Lorentzian version known as
\emph{causal dynamical triangulations} that has achieved extremely good results in the reconstruction of some of the features of general relativity (such as the four dimensionality of space-time) in the ``large scale limit'' 
(see J.~Ambj\o rn-J.~Jurkiewicz-R. Loll~\cite{AJL,AJL2} and references therein).  
\item 
\emph{Covariant canonical quantization} is a non-perturbative approach based either on the usage of field quantization via R.~Peierls brackets~\cite{Pei,BEMS} (see B.~DeWitt~\cite{DW} for details) or on covariant quantization on phase 
space~\cite{ABR}. 
\item 
\emph{Precanonical quantum gravity} is a non-perturbative covariant approach based on T.~De Donder-H.~Weyl~\cite{DeD,Wey} Hamiltonian formulation of field theory that is studied by I.~Kanatchikov (see~\cite{Kan} and references therein). 
\item 
\emph{Affine quantum gravity}, developed by J.~Klauder~\cite{Kl1,Kl3,Kl4,Kl5,Kl6,Kl7,Kl8,Kl9}, is based on a non-canonical (affine) quantization that makes heavy use of coherent states and projection operator methods~\cite{Kl0,Kl2} for dealing with quantum contraints. 
\end{itemize}
\item[b)] \emph{Relativizations of quantum mechanics.} \\
In this case we are forcing as much as possible of the formalism required by general covariance on quantum mechanics (eventually modifying it if necessary).
Although the proposal is very natural, there are almost no developed programs following this approach, probably because traditionally ``quantization'' has always been the standard route;
\begin{itemize}
\item 
Following seminal ideas by P.~Dirac~\cite{Dir}, J.-M.~Souriau~\cite{Sou} and G.~Esposito-G.~Gionti-C.~Stornaiolo~\cite{EGS}, 
C.~Rovelli~\cite{Ro12,Ro14,Ro15,Ro16,Ro} has developed a covariant formulation of classical and quantum mechanics that is appropriate for the needs of quantum relativity. 
\item
K.~Fredenhagen-R.~Haag~\cite{FH}, and more recently R.~Brunetti-K.~Fre\-den\-hagen-R.~Verch~\cite{BFV} have been studying the problem in the context of algebraic quantum field theory.
\item
A few researchers, among them B.~Mielnik~\cite{Mi} and more recently  A.~Ashtekar-T.~Schilling~\cite{AS},
C.~Brody-L.~Hugston~\cite{BH1,BH2}, have been trying to modify the usual phase-space of quantum mechanics (the K\"ahler manifold given by the projective space of a separable Hilbert space with the Fubini-Study metric) in order to allow more ``geometrical variability'' in the hope to facilitate the confrontation with general relativity.
\item
The ``consistent histories formulation'' of quantum mechanics elaborated by R.~Griffiths~\cite{Gr}, R.~Omnes~\cite{Om1,Om2}, 
M.~Gell-Mann-J.~Hartle~\cite{Hart}, and more recently the ``history projection operator theory'' developed by C~Isham-N.~Linden-N.~Savvidou-S.~Shreckenberg~\cite{I3,I4,IL0,IL1,ILSS}, provides another covariant generalization of quantum mechanics that is suitable for quantum gravity~\cite{IS1,IS2,Sav1,Sav2,Sav3,Sav4,Sav5,Sav6,Sav7}.  
\item
Some proposal to modify quantum mechanics in a ``relational'' or ``covariant way'' starting with  H.~Everett-J.~Wheeler and more recently with C.~Rovelli~\cite{Ro3,Ro12},  C.~Rovelli-M.~Smerlak~\cite{RS} or with the use of categories/topoi
(L.~Crane~\cite{Cr, Cr2},
J.~Butterfield-C.~Isham~\cite{BI3,BI4,BI5,I6,I7}, C.~Isham-A.~D\"oring~\cite{DI1,DI2,DI3,DI4,DI5,Dor1,Dor2}, C.~Flori~\cite{Fl})  in order to make it suitable for quantization of general relativity (either in the case of loop quantum gravity program of other more radical approaches) can be considered also in this category.
\end{itemize}
\item[c)] \emph{General relativity as an emergent theory.} \\
Here quantum mechanics and quantum field theory are considered as basic and general relativity is obtained as an approximation from a fundamental theory. These kind of approaches pioneered by A.~Sakharov's ``induced gravity'' have always been the most fashionable among particle physicists and are now gaining momentum also among ``relativists''.
\begin{itemize}
\item \label{it: c}
\emph{String theory} in all of its variants is the most popular approach to quantum gravity. We refer to M.~Green-J.~Schwarz-E.~Witten~\cite{GSW}, J.~Polchinski~\cite{Pol} and K.~Baker-M.~Baker-J.~Schwarz~\cite{BBS} as standard references. 
\item
\emph{Analog gravity} and other models of general relativity based on quantum solid state physics, acoustic, hydrodynamics. For a review, see for example G.~Volovik~\cite{Vo1,Vo2,Vo3} and  C.~Barcel\'o-S.~Liberati-M.~Visser~\cite{BLV}.
\item
\emph{H\v orava gravity}~\cite{Hor} is a non-relativistic quantum field theory of gravitons in $3+1$-dimensions where Lorentz invariance and relativity emerge only as approximations in the long scale limit. 
\item 
\emph{Emergent Gravity}: inspired by the partial achievements of ``analog gravity'' a new cluster area of research in gravity, seen as an emergent large-scale phenomenon, is gaining momentum (see for example the papers by F.~Girelli-S.~Liberati-L.~Sindoni~\cite{GLS,GLS2,GLS3,GLS4,GLS5}).
Some of the more recent developments of the ``path integral'' approach to quantum gravity such ``group field theory'' by 
D.~Oriti~\cite{Or,Or3,Or2} or ``causal dynamical triangulations'' by J.~Ambj\o rn-J.~Jurkiewicz-R.~Loll~\cite{AJL,AJL2} as well as 
some more radical proposals to obtain space (but not time!) and gravity as emergent from a quantum substratum such as ``internal quantum gravity'' by O.~Dreyer~\cite{D1,D2,D3,D4,D5,D6}, ``quantum causal history'' and ``quantum graphity'' by F.~Markopoulou and collaborators~\cite{Mar,Mar2,Mar3,Mar4,HMS,KM,KMS,KMSe}, ``causal sets'' by R.~Sorkin~\cite{So0,So1,So2} can also be consided in this category.  
\end{itemize}
\item[d)] \emph{Quantum mechanics as an emergent theory} (without modification of general relativity). \\
Very few people have been trying this road, probably because everyone is expecting that a classical theory (as general relativity is) should be subject to quantum  modifications in the small distances regime, there are anyway some incomplete ideas in this direction:
\begin{itemize}
\item
G.~t'Hooft~\cite{tH1,tH2,tH3,tH4,tH5,tH6} is proposing models to replace quantum mechanics with a classical fundamental deterministic theory.
\item 
The developments of geometrodynamics, as described in the recent review by D.~Giulini~\cite{Giu}, suggest the possibility to recover at least some of the properties of matter from pure geometry.  
\item 
The theory of geons (H.~Hadley~\cite{Ha1,Ha2,Ha3}), tries to simulate the quantum behaviour of elementary particles starting with localized geometrical structures on the Lorentzian manifolds of general relativity.
\item
L.~Smolin~\cite{Sm2} has recently considered the possibility that quantum mechanics might arise as a stochastic theory induced by non-local variables.
\item 
E.~Prugovecki~\cite{Pr1,Pr2,Pr3} also proposed an approach to quantum mechanics through stochastic processes in a general relativistic geometrical setting. 
\item
The theory of gravitational induced collapse of the quantum wave function by R.~Penrose (see~\cite{Pe} and references therein) can be considered in this category. 
\end{itemize}
\item[e)] \emph{Pregeometrical approaches} (suggested by J.~Wheeler)
are alternative approaches that require at least some basic modifications of general relativity and quantum mechanics that might both ``emerge'' by some deeper dynamic of degrees of freedom not necessarily related to any macroscopic geometrical entity.
Most of these theories are at least partially background-independent (depending on the amount of ``residual'' geometrical structure used to define their kinematic). The main problems arising in pregeometrical theories is usually the description of an appropriate dynamic and the recovery from it of some ``approximate'' description of general relativity and ordinary quantum physics in the ``macroscopic'' limit.
The proposals that can be listed in this category are extremely heterogeneous and they might range from ``generalizations'' of other more conservative approaches:
\begin{itemize}
\item
algebraic quantum gravity: a generalization of loop quantum gravity recently developed by K.~Giesel-T.~Thiemann~\cite{GT1,GT2,GT3,GT4},
\item
group field theory quantum gravity: a powerful extension of the path integral approach to quantum gravity proposed by 
D.~Oriti~\cite{Or,Or2},
\end{itemize}
to more radical paths (that we collect here just for the benefit of the interested reader):
\begin{itemize}
\item
twistor theory (R.~Penrose~\cite{PeR, Pe}),
\item
quantum code (D.~Finkelstein~\cite{Fi1,Fi2}),
\item
causal sets (R.~Sorkin~\cite{So0,So1,So2}),
\item
causaloids (L.~Hardy~\cite{Har1,Har2,Har3,Har4,Har5,Har6}),
\item
computational approach (S.~Lloyd~\cite{Ll1,Ll2}),
\item
internal quantum gravity (O.~Dreyer~\cite{D1,D2,D3,D4,D5,D6}),
\item
quantum causal history
(F.~Markopoulou~\cite{Mar,Mar2,Mar3,Mar4}, E.~Hawkins-F.~Markopoulou-H.~Sahlman~\cite{HMS}, 
D.~Kribs-F.~Makopoulou~\cite{KM});  
quantum graphity (T.~Konopka-F.~Markopoulou-L.~Smolin-S.~Severini~\cite{KMS,KMSe}),
\item
abstract differential geometry (A.~Mallios~\cite{Mal1,Mal2,Mal3}, J.~Raptis~\cite{Ra1,Ra2,Ra3,Ra4,Ra5,Ra6}, A.~Mallios-J.~Raptis~\cite{MR1,MR2,MR3,MR4}, A.~Mallios-E.~Rosinger~\cite{MRo}),
\item
categorical approaches
(J.~Baez~\cite{Ba2,Ba3,Ba4}, 
L.~Crane~\cite{Cr,Cr2,Cr3,Cr4}, J.~Morton~\cite{Mor3}, 
J.~Butterfield-C.~Isham~\cite{BI3,BI4,BI5}, C.~Isham~\cite{I7,I18,I19,I20}, 
A.~D\"oring-C.~Isham~\cite{DI1,DI2,DI3,DI4,DI5,Dor1,Dor2}), 
C.~Flori~\cite{Fl}, 
\item
non-commutative geometry approaches:\footnote{See also the recent papers by B.~Booss-Bavnbek-G.~Esposito-M.~Lesch~\cite{BEL} and F.~M\"uller-Hoissen~\cite{M-H} for more detailed and alternative surveys on noncommutative geometry in gravity.}
\begin{itemize}
\item
via derivations on non-commutative (groupoid) algebras: J.~Madore~\cite{Ma},
M.~Heller-Z.~Odrzygozdz-L.~Pysiak-W.~Sasin~\cite{HOS,HS1,HS2,HS3,HS4,HS5,HS6,HPS1,HPS2,HPS3,HPS4,HOPS1,HOPS2,HOPS3,HOPS4},
\item
via deformation quantization (Moyal-Weyl):
P.~Aschieri and collaborators~\cite{As1,As2,ADMW,ADMSW, ABDMSW},
\item
via quantum groups:
S.~Majid~\cite{Maj3,Maj6,Maj4,Maj5},
\item
via A.~Connes' non-commutative geometry:   \\
M.~Paschke~\cite{P}, A.~Connes-M.~Marcolli~\cite{CM6}.
\end{itemize}
\end{itemize}
\end{itemize}
Since we are here mainly interested in A.~Connes' non-commutative geometry, we are going to conclude by examining a bit more in detail the situation as regards its possible applications to quantum gravity.

\subsubsection{A.~Connes' Non-commutative Geometry and Gravity}\label{sec: qgrav}

It is often claimed that non-commutative geometry will be a key ingredient (a kind of quantum version of Riemannian geometry) for the formulation of a fundamental theory of quantum gravity (see for example L.~Smolin~\cite{Sm} and P.~Martinetti~\cite{Mart}) and actually non-commutative geometry is often listed among the current alternative approaches to quantum gravity.

In reality, with the only notable exceptions of the extremely interesting programs outlined in M.~Paschke~\cite{P} and in A.~Connes-M.~Marcolli~\cite{CM6}, a foundational approach to quantum physics based on A.~Connes' non-commutative geometry has never been proposed.
So far, most of the current applications of A.~Connes' non-commutative geometry to (quantum) gravity have been limited to:
\begin{itemize}
\item
the study of some ``quantized'' example: C.~Rovelli~\cite{Ro8}, F,~Besnard~\cite{Be},
\item
the use of its mathematical framework for the reformulation of classical (Euclidean) general relativity: D.~Kastler~\cite{K2}, A.~Chamseddine-G.~Felder-J.~Fr\"ohlich~\cite{CFF}, \\  W.~Kalau-M.~Walze~\cite{KW}, C.~Rovelli-G.~Landi~\cite{LR1,LR2,Lan2},
\item
attempts to use its mathematical framework ``inside'' some already established theories such as string theory (A.~Connes-M.~Douglas-A.~Schwarz~\cite{CDS}, J.~Fr\"ohlich-O.~Grandjean-A.~Recknagel~\cite{FGR3} and J.~Brodzki-V.~Mathai-J.~Rosenberg-R.~Szabo~\cite{BMRS}) and loop gravity 
(J.~Aastrup-J.~Grimstrup-R.~Nest-M.~Paschke~\cite{AG1,AG2,AGN1,AGN2,AGN3,AGN4,AGNP}, 
F.~Girelli-E.~Livine~\cite{GL}),
\item
the formulation of Hamiltonian theories of gravity on globally hyperbolic cases, where only the ``spacial-slides'' are described by non-commutative geometries: \\ E.~Hawkins~\cite{Haw}, T.~Kopf-M.~Paschke~\cite{KP1,KP2,Ko}.
\end{itemize}

\label{sec: modular}

In a slightly different direction, there are some important areas of research that are somehow connected to the problems of quantum gravity and that seem to suggest a more prominent role of Tomita-Takesaki modular theory\footnote{The original ideas about modular theory were developed by M.~Tomita~\cite{To1,To2}. We refer to the texts by M.~Takesaki~\cite{T}, B.~Blackadar~\cite{Bl} for a modern mathematical introduction and to O.~Bratteli-D.~Robinson~\cite{BR}, R.~Haag~\cite{H} for a more physics oriented presentation. Excellent updated reviews on the relevance of modular theory in quantum physics are given by 
S.~J.~Summers~\cite{Su}, H.-J.~Borchers~\cite{Bo} and the recent works by D.~Guido~\cite{G} and F.~Lled\'o~\cite{Lle2} 
 (but see also R.~Longo \cite{L0}). Outside the realm of operator algebras, Tomita-Takesaki theorem for classical statistical mechanical systems has been discussed by G.~Gallavotti-M.~Pulvirenti~\cite{GP} and a strictly related  correspondence between modular theory and Poisson geometry has been pointed out by A.~Weinstein~\cite{W}.}
in quantum physics (and in particular in the physics of gravity):
\begin{itemize}
\item
Since the work of W.~Unruh~\cite{U}, it has been conjectured the existence of a deep connection between gravity (equivalence principle), thermal physics (hence Tomita-Takesaki and KMS-states) and quantum field theory; this idea has not been fully exploited so far. This line of thought is actually reinforced by the works on thermodynamical derivation of Einstein equation by T.~Jacobson~\cite{Jac} (see also R.~Brustein-M.~Hadad~\cite{BH} and M.~Parikh-S.~Sarkar~\cite{PaS}). 
\item
Starting from the works of J.~Bisognano-E.~Wichmann~\cite{BW1,BW2}, G.~Sewell~\cite{Sew} and more recently, H.~J.~Borchers~\cite{Bo1}, there is mounting evidence that Tomita-Takesaki modular theory should  play a fundamental role in the ``spectral  reconstruction'' of the space-time information from the algebraic setting of states and observables.
The most interesting results in this direction have been obtained so far:
\item[-]
in the theory of ``half-sided modular inclusions'' and modular intersections (see H.-J.~Borchers~\cite{Bo} and references therein, H.~Araki-L.~Zsido~\cite{ArZs});
\item[-]
in the ``geometric modular action'' program (see for more details D.~Buchholz-S.~J.~Summers~\cite{BuS,BS2}, D.~Buchholz-M.~Florig-S.~J.~Summers~\cite{BFS}, D.~Buchholz-O.~Dreyer-M.~Florig-S.~J.~Summers~\cite{BDFS}, 
S.~J.~Summers-R.~White~\cite{SuW}); 
\item[-]
in ``modular nuclearity'' (see for more details R.~Haag~\cite{H} and, for recent applications to the ``form factor program'', D.~Buchholz-G.~Lechner~\cite{BL,Le1,Le2,Le3,Le4,Le5}, D.~Buchholz-S.~J.~Summers~\cite{BS3}); 
\item[-]
in the ``modular localization program'' (see B.~Schroer-H.-W.~Wiesbrock~\cite{Sc1,Sc2,SW1,SW2}, R.~Brunetti-D.~Guido-R.~Longo~\cite{BGL}, F.~Lled\'o~\cite{Lle}, J.~Mund-B.~Schroer-J.~Yngvanson~\cite{MSY} and N.~Pinamonti~\cite{Pi}).
\item
Starting with the construction of cyclic cocycles from supersymmetric quantum field theories by A.~Jaffe-A.~Lesniewski-K.~Osterwalder~\cite{JLO,JLO2}, there has always been a constant interest in the possible deep structural relationship between supersymmetry, modular theory of type III von Neumann algebras and non-commutative geometry (see D.~Kastler~\cite{K2} and A.~Jaffe-O.~Stoytchev~\cite{J,JS}). Some deep results by R.~Longo~\cite{L} established a bridge between the theory of superselections sectors and cyclic cocycles obtained by super-KMS states. 
The recent work by D.~Buchholz-H.~Grundling~\cite{BG,BGr2} opens finally a way to construct super-KMS functionals and spectral triples in algebraic quantum field theory (see S.~Carpi-R.~Hillier-Y.~Kawahigashi-R.~Longo~\cite{CHKL}).
\item
In the context of C.~Rovelli ``thermal time hypothesis''~\cite{Ro} in quantum gravity,  A.~Connes-C.~Rovelli~\cite{CR} (see also P.~Martinetti-C.~Rovelli~\cite{MR} and P.~Martinetti~\cite{Mart0,Mart3}) have been using Tomita-Takesaki modular theory in order to induce a macroscopic time evolution for a relativistic quantum system.
\item
A.~Connes-M.~Marcolli~\cite{CM6} with the ``cooling procedure'' are proposing to examine the operator algebra of observables of a quantum gravitational system, via modular theory, at ``different temperatures'' in order to extract by ``symmetry breaking'' an emerging geometry.
\end{itemize}

\begin{itemize}
\item[$\looparrowright$]
The idea that space-time might be spectrally reconstructed, via non-commutative geometry, from Tomita-Takesaki modular theory applied to the algebra of physical observables was elaborated in 1995 by one of the authors (P.B.) and independently (motivated by the possibility to obtain cyclic cocycles in algebraic quantum field theory from modular theory) by R.~Longo~\cite{L1}. Since then this conjecture is still the main subject and motivation of our investigation~\cite{B,BCL}.
\end{itemize}
Similar speculations on the interplay between modular theory and (some aspects of) space-time geometry have been suggested by S.~Lord~\cite[Section~VII.3]{Lo} and by M.~Paschke-R.~Verch~\cite[Section~6]{PV}.
\begin{itemize}
\item[$\looparrowright$]
One of the authors (R.C.) has raised the somehow puzzling question whether it is possible
to reinterpret the one parameter group of modular automorphisms as a renormalization (semi-)group in physics.
The connection with P.~Cartier's idea of a ``universal Galois group''~\cite{Car}, currently developed by A.~Connes-M.~Marcolli, is extremely intriguing.
\end{itemize}

\subsubsection{A Proposal for (Modular) Algebraic Quantum Gravity.}\label{sec: AQG}

Our ongoing research project~\cite{B1,BCL,B3}~\footnote{P.~Bertozzini, R.~Conti, W.~Lewkeeratiyutkul, Modular Algebraic Quantum Gravity, work in progress.}
is aiming at the construction of an \emph{algebraic theory of quantum gravity} in which ``non-commutative'' space-time is spectrally reconstructed from Tomita-Takesaki modular theory.

What we propose is to develop an approach to the foundations of quantum physics technically based on algebraic quantum theory (operator algebras) and A.~Connes' non-commutative geometry. The research is building on the experience already gained in our previous/current mathematics research plans on ``modular spectral triples in non-commutative geometry and physics''~\cite{BCL}\footnote{Partially supported by the Thai Research Fund~\hlink{http://www.trf.or.th/research/project.asp?projectid=RSA4580030}{TRF project RSA4580030}.} and on ``categorical non-commutative geometry'' and
is conducted in the standard of mathematical rigour typical of the tradition of mathematical physics' research in algebraic quantum field theory~\cite{A,H}.

In the mathematical framework of A.~Connes' non-commutative geometry, we are addressing the problem of the ``spectral reconstruction'' of ``geometries'' from the underlying operational data defined by ``states'' over ``observables' \cs-algebras'' of physical systems.
More specifically:
\begin{itemize}
\item[$\looparrowright$] 
Building on our previous research on ``modular spectral triples''\footnote{P.~Bertozzini, R.~Conti, W.~Lewkeeratiyutkul, Modular Spectral Triples, in preparation.} and on recent results on semi-finite spectral triples developed by A.~Carey-J.~Phillips-A.~Rennie-F.~Sukhochev~\cite{CPR,CPR2,CPR3,CPR4}\footnote{See also M.~Laca--S.~Neshveyev~\cite{LN} and A.~Carey-S.~Neshveyev-R.~Nest-A.~Rennie~\cite{CNNR}.},
we make use of Tomita-Takesaki modular theory of operator algebras to associate, to suitable states $\omega$ over involutive normed algebras $\As$, non-commutative geometrical objects $(\As_\omega,\Hs_\omega,D_\omega)$ that are only formally similar to A.~Connes' spectral-triples and where the ``Dirac operator'' $D_\omega$, that is usually taken as the modular Hamiltonian 
$K_\omega=\log \Delta_\omega$, satisfies the modular invariance property $\Delta^{it}D_\omega\Delta^{-it}=D_\omega$  
(for some more details see the slides~\cite[Pages~74-77]{B3}). 
\item[$\looparrowright$]
We are now developing\footnote{P.~Bertozzini, Algebraic Formalism for Rovelli Quantum Theory,
in preparation.} an ``event'' interpretation of the formalism of states and observables in algebraic quantum physics that is in line with C.~Isham's ``history projection operator theory''~\cite{I3,I4,IL0,IL1,ILSS} and C.~Rovelli's ``relational/relativistic quantum mechanics''~\cite{Ro3} (for additional details see the slides~\cite[Pages~78-81]{B3}). 
\item[$\looparrowright$]
Making contact with our current research project on ``categorical non-commutative geometry'' and with other projects in categorical quantum gravity (J.~Baez~\cite{Ba3,Ba4} and L.~Crane~\cite{Cr,Cr2}), we plan to generalize the diffeomorphism covariance group of general relativity in a categorical context and use it to ``identify'' the degrees of freedom related to the spatio-temporal structure of the physical system (more details can be found in the slides~\cite[Pages~82-84]{B3}). 
\item[$\looparrowright$]
Techniques from ``decoherence/einselection'' (H.~Zeh~\cite{Ze}, W.~Zurek~\cite{Zu}),
``emergence/noiseless subsystems'' (O.~Dreyer~\cite{D1,D2,D3}, F.~Markopoulou~\cite{Mar,Mar2,Mar3,KoM}), superselection (I.~Ojima~\cite{O,O2,OT}) and the ``cooling'' procedure developed by A.~Connes-M.~Marcolli~\cite{CM6} are expected to be relevant in order to extract from our spectrally defined non-commutative geometries, a macroscopic space-time for the pair state/system and its ``classical residue''.
\item[$\looparrowright$]
Possible reproduction of quantum geometries already defined in the context of loop quantum gravity (T.~Thiemann~\cite{Th1,Th4} and J.~Aastrup-J.~Grimstrup-R.~Nest-M.~Paschke~\cite{AG1,AG2,AGN1,AGN2,AGN3,AGN4,AGNP}) or in S.~Doplicher-J.~Roberts-K.~Fredenhagen models~\cite{DFR0,DFR,Do1,Do2,Do} will be investigated.
\end{itemize}
If partially successful, the project will have a significant fallout:
a background-independent powerful approach to ``quantum relativity'' that is suitable for the purpose of unification of physics, geometry and information theory that lies ahead.

\appendix 

\section*{Appendix: Some Recent Developments} \label{appendix}

The first version of this arXiv preprint was written in November 2007 and this second corrected and expanded version was actually prepared in November 2009. Now, in December 2011, after more than two years, a few notable developments occurred, but we decided, for this final arXiv version, to ``freeze'' the bibliographical references directly discussed in the paper to October 2009, limiting our revision of the main text to correction of misprints and updates only of those bibliographical sources already appeared in preprint before November 2009. 

As a partial remedy, in this appendix we provide, for the interested reader, (a very selective choice of) a few additional bibliographical references to recently appeared works mainly related to A.~Connes' spectral triples in non-commutative differential geometry. 

\medskip 

On the ``Riemannian version of spectral triples'': 
\begin{itemize}
\item 
Lord, Steven; Rennie Adam; Varilly, Joseph C., 
Riemannian Manifolds in Noncommutative Geometry, 
\hlink{http://arxiv.org/abs/1109.2196}{arXiv:1109.2196v1}.
\end{itemize}
For locally compact spectral triples: 
\begin{itemize}
\item
Carey, Alan; Gayral, Victor; Rennie, Adam; Sukochev, Fedor, 
Index Theory for Locally Compact Noncommutative Geometries
\hlink{http://arxiv.org/abs/1107.0805}{arXiv:1107.0805v1}.
\end{itemize}

On Lorentzian non-commutative geometry: 
\begin{itemize}
\item
Verch Rainer, 
Quantum Dirac Field on Moyal-Minkowski Spacetime - Illustrating Quantum Field Theory over Lorentzian Spectral Geometry, 
\hlink{http://arxiv.org/abs/1106.1138}{arXiv:1106.1138v1}.
\item
Franco, Nicolas, 
Lorentzian Approach to Noncommutative Geometry, \\ 
\hlink{http://arxiv.org/abs/1108.0592}{arXiv:1108.0592v1}.
\end{itemize}

Variations of Connes' reconstructions theorem for almost commutative spectral triples: 
\begin{itemize}
\item 
\'Ca\'ci\'c, Branimir, 
A Reconstruction Theorem for Almost-commutative Spectral \\ Triples, 
\hlink{http://arxiv.org/abs/1101.5908}{arXiv:1101.5908v3}.
\end{itemize}

On spectral characterization of isometries: 
\begin{itemize}
\item
Cornelissen, Gunther; de Jong, Jan Willem, 
The Spectral Length of a Map Between Riemannian Manifolds, 
\hlink{http://arxiv.org/abs/1007.0907}{arXiv:1007.0907v3}. 
\end{itemize}

Spectral triples on crossed products: 
\begin{itemize}
\item
Bellissard, Jean; Marcolli, Matilde; Reihani, Kamran, 
Dynamical Systems on Spectral Metric Spaces, 
\hlink{http://arxiv.org/abs/1008.4617}{arXiv:1008.4617v1}.
\item
Hawkins, Andrew; Skalski, Adam; White, Stuart; Zacharias, Joachim, 
Spectral Triples on Crossed Products Arising from Equicontinuous Actions, \\ 
\hlink{http://arxiv.org/abs/1103.6199}{arXiv:1103.6199v3}.
\end{itemize} 

Further works on application of non-commutative geometry to the standard model and physics are: 
\begin{itemize}
\item 
Chamseddine, Ali; Connes Alain (2010). 
Noncommutative Geometry as a Framework for Unification of all Fundamental Interactions including Gravity. Part I
\textit{Fortsch.~Phys.} 58, 553-600,  
\hlink{http://arxiv.org/abs/1004.0464}{arXiv:1004.0464v1}. 
\item 
Chamseddine, Ali; Connes Alain, 
Space-Time from the Spectral Point of View, 
\hlink{http://arxiv.org/abs/1008.0985}{arXiv:1008.0985v1}. 
\item
Chamseddine, Ali; Connes Alain (2011). 
Noncommutative Geometric Spaces with Boundary: Spectral Action,  
\textit{J.~Geom.~Phys.} 61 n.~1, 317-332,
\hlink{http://arxiv.org/abs/1008.3980}{arXiv:1008.3980v1}.
\item 
van den Dungen, Koen; van Suijlekom Walter,  
Electrodynamics from Noncommutative Geometry
\hlink{http://arxiv.org/abs/1103.2928}{arXiv:1103.2928v1}. 
\item 
Boeijink, Jord; van Suijlekom, Walter, 
The Noncommutative Geometry of Yang-Mills Fields, 
\hlink{http://arxiv.org/abs/1008.5101}{arXiv:1008.5101v1}.
\item
van den Broek, Thijs; van Suijlekom, Walter, 
Supersymmetric QCD and Noncommutative Geometry
\hlink{http://arxiv.org/abs/1003.3788}{arXiv:1003.3788v1}. 
\item
Bhowmick, Jyotishman; D'Andrea, Francesco; Das, Biswarup; Dabrowski, Ludwik, 
Quantum Gauge Symmetries in Noncommutative Geometry, 
\hlink{http://arxiv.org/abs/1112.3622v1}{arXiv:1112.3622v1}.
\end{itemize}
For the study of Connes' spectral distance see the following papers and references therein: 
\begin{itemize}
\item 
Cagnache, Eric; D'Andrea, Francesco; Martinetti, Pierre; Wallet, Jean-Christophe (2011).  
The Spectral Distance on the Moyal Plane,
\textit{J.~Geom.~Phys.} 61, 1881-1897, 
\hlink{http://arxiv.org/abs/0912.0906}{arXiv:0912.0906v3}.
\item 
Martinetti, Pierre; Mercati, Flavio; Tomassini, Luca, 
Minimal Length in Quantum Space and Integrations of the Line Element in Noncommutative Geometry
\hlink{http://arxiv.org/abs/1106.0261}{arXiv:1106.0261v1}.
\item 
Martinetti, Pierre; Tomassini, Luca, 
Noncommutative Geometry of the Moyal Plane: Translation Isometries and Spectral Distance Between Coherent States, 
\\ 
\hlink{http://arxiv.org/abs/1110.6164}{arXiv:1110.6164v1}.
\end{itemize}

Semi-finite and modular spectral triples are treated in: 
\begin{itemize}
\item 
Carey, Alan; Phillips, John; Putnam, Ian; Rennie Adam (2011). 
Families of Type III KMS States on a Class of C*-algebras Containing $O_n$ and $\mathcal{Q}_\N$, 
\textit{J.~Funct.~Anal.} 260 n.~6, 1637-1681, 
\hlink{http://arxiv.org/abs/1001.0424}{arXiv:1001.0424v1}.
\item
Lai, Alan, 
On Type II Noncommutative Geometry and the JLO Character, \\ 
\hlink{http://arxiv.org/abs/1003.4226}{arXiv:1003.4226v1}.
\item
Rennie, Adam; Senior, Roger, 
The Resolvent Cocycle in Twisted Cyclic Cohomology and a Local Index Formula for the Podles Sphere, 
\hlink{http://arxiv.org/abs/1111.5862}{arXiv:1111.5862v1}.
\item
Rennie, Adam; Sitarz, Andrzej; Yamashita, Makoto, 
Twisted Cyclic Cohomology and Modular Fredholm Modules, 
\hlink{http://arxiv.org/abs/1111.6328}{arXiv:1111.6328v1}.
\item 
Kaad, Jens, 
On Modular Semifinite Index Theory, 
\hlink{http://arxiv.org/abs/1111.6546}{arXiv:1111.6546v1}.
\end{itemize} 

On ``Morita morphisms'' of spectral triples, beside Bram Mesland work (now already cited in the main paper): 
\begin{itemize} 
\item 
Kaad, Jens; Lesch, Matthias,   
Spectral Flow and the Unbounded Kasparov Product,
\hlink{http://arxiv.org/abs/1110.1472}{arXiv:1110.1472v1}.
\end{itemize}

As regards non-commutative geometrical approaches to (loop) quantum gravity: 
\begin{itemize}
\item
Denicola, Domenic; Marcolli, Matilde; Zainy al-Yasry, Ahmad (2010).  
Spin Foams and Noncommutative Geometry. 
\textit{Classical Quantum Gravity} 27 n.~20, 205025, 53 pp. 
\hlink{http://arxiv.org/abs/1005.1057}{arXiv:1005.1057v1}.
\item
Gracia-Bondia Jose, 
Notes on "Quantum Gravity" and Non-commutative Geometry, 
\hlink{http://arxiv.org/abs/1005.1174}{arXiv:1005.1174v1}.
\item
Lai, Alan, 
The JLO Character for The Noncommutative Space of Connections of Aastrup-Grimstrup-Nest
\hlink{http://arxiv.org/abs/1010.5226}{arXiv:1010.5226v1}.
\item 
Aastrup, Johannes; Grimstrup, Jesper M\o ller; Paschke, Mario (2011).  
Quantum Gravity Coupled to Matter via Noncommutative Geometry, \\ 
\textit{Classical Quantum Gravity} 28 n.~7, 075014, 10 pp., 
\hlink{http://arxiv.org/abs/1012.0713}{arXiv:1012.0713v1}.
\item 
Aastrup, Johannes; Grimstrup, Jesper M\o ller; 
From Quantum Gravity to Quantum Field Theory via Noncommutative Geometry
\hlink{http://arxiv.org/abs/1105.0194}{arXiv:1105.0194v1}. 
\end{itemize}

Our work on modular algebraic quantum gravity has received a more detailed treatment in the paper:
\begin{itemize}
\item
Bertozzini, Paolo; Conti, Roberto; Lewkeeratiyutkul, Wicharn (2010). 
Modular Theory, Non-commutative Geometry and Quantum Gravity. 
\textit{SIGMA Symmetry Integrability Geom.~Methods Appl.} 6 paper 067, 47 pp.
\hlink{http://arxiv.org/abs/1007.4094}{arXiv:1007.4094v2}. 
\end{itemize}

For studies recently appeared on the usage of Tomita-Takesaki modular theory in quantum physics and loop quantum gravity, that although not directly related to spectral triples in non-commutative geometry, might have deep impact on our approach to modular algebraic quantum gravity see the following preprints and the references therein: 
\begin{itemize}
\item 
Asselmeyer-Maluga, Torsten; Krol, Jerzy, 
Constructing a Quantum Field Theory from Spacetime, 
\hlink{http://arxiv.org/abs/1107.3458}{arXiv:1107.3458v1}. 
\item 
Kaminski, Diana, 
Algebras of Quantum Variables for Loop Quantum Gravity, I. Overview, 
\hlink{http://arxiv.org/abs/1108.4577}{arXiv:1108.4577v1}.
\item 
Kostecki, Ryszard, 
Information Dynamics and New Geometric Foundations of Quantum Theory, 
\hlink{http://arxiv.org/abs/1110.4492}{arXiv:1110.4492v3}. 
\end{itemize}

\end{document}